\documentclass[11pt]{amsart} 

\usepackage[latin9]{inputenc}
\usepackage{amssymb}
\usepackage{latexsym}
\usepackage{amsmath}
\usepackage{amsthm}
\usepackage{amsfonts}

\usepackage{bbm}
\allowdisplaybreaks

\newtheorem{theorem}{Theorem}[section]

\newtheorem{lemma}[theorem]{Lemma}
\newtheorem*{lemma*}{Lemma}
\newtheorem{corollary}[theorem]{Corollary}
\newtheorem{proposition}[theorem]{Proposition}
\newtheorem*{proposition*}{Proposition}

\theoremstyle{definition}

\newtheorem*{acknowledgements*}{Acknowledgements}
\theoremstyle{remark}
\newtheorem{remark}[theorem]{Remark}

\numberwithin{equation}{section}

\usepackage{color}

\newcommand{\R}{\mathbb R}
\newcommand{\N}{\mathbb N}
\newcommand{\C}{\mathbb C}
\newcommand{\T}{\mathbb T}
\newcommand{\D}{\mathbb D}
\newcommand{\Z}{\mathbb Z}
\newcommand{\ct}{\mathcal{C}_t}
\newcommand{\cno}{\C^{\N_0}}
\newcommand{\no}{\N_0}
\newcommand{\ce}{\mathcal{C}}

\newcommand{\ces}{Ces\`{a}ro }
\newcommand{\subi}{{n=0}^\infty}
\newcommand{\txop}{{X\to X}}
\newcommand{\ee}{\mathcal{E}}
\newcommand{\np}{{N^p(\Z)}}
\newcommand{\xo}{x^{(\text{o})}}
\newcommand{\xe}{x^{(\text{e})}}

\title
[Generalized Ces\`{a}ro operators]
{Fine spectra and compactness of   
\\ generalized Ces\`{a}ro operators \\ in  
Banach lattices in $\cno$}

\author[G.~P.  Curbera]{Guillermo P. Curbera}
\address{Facultad de Matem\'aticas \& IMUS,
Universidad de Sevilla, 
Calle Tarfia s/n,  Sevilla 41012, Spain}
\email{curbera@us.es}

\author[W.~J. Ricker]{Werner J. Ricker}
\address{Math.--Geogr. Fakult\"at, Katholische Universit\"at
Eichst\"att--Ingolstadt, D--85072 Eichst\"att, Germany}
\email{werner.ricker@ku.de}

\thanks{The first author acknowledges the support  of 
PGC2018-096504-B-C31, FQM-262 and Feder-US-1254600 (Spain).}

\date{\today}

\begin{document}

\begin{abstract}
The  generalized Ces\`{a}ro operators $\ct$, for $t\in[0,1)$,
introduced  in the 1980's by  Rhaly, are natural analogues 
of the classical Ces\`{a}ro averaging operator $\ce_1$ and act in various 
Banach sequence spaces  $X\subseteq \cno$.  In this paper we concentrate
on a certain class of   Banach lattices for the coordinate-wise order,
which includes  all  separable, rearrangement invariant sequence spaces, 
various weighted $c_0$ and $\ell^p$ spaces 
 and many others. In such Banach lattices $X$
the operators $\ct$, for $t\in[0,1)$, are always compact (unlike $\ce_1$)
and a full description of their point, continuous and residual spectrum is given.
Estimates for the operator norm of $\ct$ are also presented.
\end{abstract} 

\maketitle


\section{Introduction}
\label{S1}


The  \textit{(discrete) generalized \ces operators} 
$\ct$, for $t\in(0,1)$, were first investigated by Rhaly, \cite{Rh}. The action of $\ct$
from $\cno$ into itself (where $\no:=\{0,1,2\dots\}$) is given by
\begin{equation}\label{11}
\ct x:=\left(\frac{t^nx_0+t^{n-1}x_1+\cdots+x_n}{n+1}\right)_{n=0}^\infty,
\end{equation}
for each $x=(x_n)_{n=0}^\infty\in\cno$.  Note that if $t=0$, then $\mathcal{C}_0$ is the diagonal 
operator
\begin{equation}\label{12}
D x:=\left(\frac{x_n}{n+1}\right)_{n=0}^\infty,
\quad x \in\cno,
\end{equation}
and for $t=1$ the operator $\mathcal{C}_1$ is the classical Ces\`{a}ro
averaging operator given by
\begin{equation}\label{13}
\mathcal{C}_1 x:=\left(\frac{x_0+x_1+\cdots+x_n}{n+1}\right)_{n=0}^\infty,
\quad x \in\cno.
\end{equation}

The spectrum and fine spectra of $\mathcal{C}_1$ have been 
extensively studied in various spaces. As a sample, 
we mention  $\ell^p$ ($1<p<\infty$),
\cite{BHS, curbera-ricker1, G, Le}, $c_0$ \cite{AB, Le,Re}, $c$  
\cite{Le}, $\ell^\infty$ \cite{P, Re},
the Bachelis spaces $N^p$ ($1<p<\infty$), \cite{curbera-ricker5}, 
$bv_0$ and $bv$ \cite{O1,O2},
weighted $\ell^p$ spaces \cite{ABR1,ABR2}, the discrete Ces\`{a}ro  spaces 
$ces_p$ ($p\in\{0\}\cup(1,\infty)$), \cite{curbera-ricker4}, 
and their dual spaces $d_s$ ($1<s<\infty$), \cite{bonet-ricker}; see also the references
in these papers. Accordingly, we will concentrate  mainly on the generalized
Ces\`{a}ro  operators $\ct$ for $t\in[0,1)$.

In the case $t\in[0,1)$, for all known cases of concrete sequence spaces 
$X\subseteq\cno$ the operator $\ct$ is compact and its fine spectra (i.e., its point, continuous and residual spectrum) are known and 
exhibit a stability independent of $X$. 
The first  results available concern $\ct$ in $X=\ell^2$ and go back to
Rhaly, \cite{Rh, Rha}.  A similar investigation occurs for  $\ell^p$ 
($1<p<\infty$) in \cite{Y-D} and for $c$ and $c_0$ 
in \cite{S-El,Y-M-D}. The recent paper \cite{S-El} also treats the operators
$\ct$ when they act in  $bv_0$, $bv$, $c$,
$\ell^1$, $\ell^\infty$, and   the Hahn space $h$; see also \cite[\S4]{MRT} for $h$.
Concerning the operator norm of $\ct$, see also
\cite{Rho1, Rho2} and the references therein.


The following conjecture, aimed at suggesting a general result which unifies and explains 
the behavior  of  $\ct$ in a large number of
different  sequence spaces, was posed in  \cite[Section 9]{S-El}.  
For the notation see Section \ref{S2}. Let 
$$
\Lambda:=\{1/(n+1):n\in\no\}.
$$

\textbf{Conjecture}:  \textit{Let $t\in[0,1)$ and $X\subseteq \cno$ be 
a  Banach sequence space containing
\begin{equation}\label{31}
\mathcal{S}:=\Big\{x\in\cno:\beta(x):=\lim_{n\to\infty}\Big|\frac{x_{n+1}}{x_n}\Big|<1\Big\}.\end{equation}
\begin{itemize}
\item[(i)] If $\ct\in\mathcal{L}(X)$, then $\sigma_{\text{p}}(\ct;X)=\Lambda$. 
Hence,  $\|\ct\|\ge1$.
\item[(ii)]  If $\ct $ is compact, then $\sigma(\ct;X)=\Lambda\cup\{0\}$.
\item[(iii)]  If $\ct $ is compact and $\mathcal{R}(\ct)$ is  dense in $X$, then 
$$
\sigma_{\text{c}}(\ct;X)=\{0\};\quad \sigma_{\text{r}}(\ct;X)=\emptyset.
$$
\item[(iv)]  If $\ct $ is compact and $\mathcal{R}(\ct)$ is not dense in $X$, then 
$$
\sigma_{\text{c}}(\ct;X)=\emptyset;\quad \sigma_{\text{r}}(\ct;X)=\{0\}.
$$
\end{itemize}
}


The aim of the present paper is to study the compactness,  the  spectrum and the fine 
spectra of the operators $\ct$, for $t\in[0,1)$.
It is shown in Theorem \ref{t-35} that  the above Conjecture  is valid 
for a large class of sequence spaces $X\subseteq\cno$. This result applies
to  \textit{solid Banach lattices} $X$ (for the coordinate-wise order) 
satisfying relatively mild requirements, mainly  on  the canonical vectors 
$\mathcal{E}=\{e_n:n\in\no\}\subseteq X$ and on the behavior 
of  the restriction to $X$ of the right-shift operator $S\colon\cno\to\cno$. 
It turns out  in such cases that  the  operators $\ct$ map $X$ 
continuously into itself, are
compact and their spectrum and fine spectra are precisely 
as stated in the above Conjecture.
Estimates for the operator norm of $\ct$ are also given.


Section \ref{S2} contains various preliminary facts on sequence spaces, in particular,
on  the extensive class of all solid Banach lattices $X\subseteq \cno$.

In Section \ref{S3} we study the action of the operators $\ct$  on  solid Banach lattices  and establish a rather general result, namely Theorem \ref{t-35}, which provides
a positive solution to the above Conjecture 
and is applicable to an extensive class of sequence
spaces. In particular, we apply Theorem \ref{t-35} to the operators $\ct$  acting on 
all separable, rearrangement invariant spaces over $\no$.

Section \ref{S6} is devoted to an investigation 
of the operators $\ct$ when they act on  the Ces\`{a}ro sequence spaces
$X=ces_p$, for $p\in\{0\}\cup(1,\infty)$, and in their dual 
spaces $d_p$, for $1\le p<\infty$. 
This family of spaces is not rearrangement invariant.
The  non-separable Ces\`{a}ro space $ces_\infty$ is also considered.

Section \ref{S8} treats the operators $\ct$ in the 
class of \textit{weighted sequence spaces}
$\ell^p(w)$, for $p\in[1,\infty)$, and $c_0(w)$,
where $w=(w_n)_\subi$ satisfies $0<w_n\downarrow0$.

In Section  \ref{S9} we study  the operators $\ct$ acting in the Bachelis spaces $X=N^p$ for $p\in(1,\infty)$, a class of spaces arising in classical 
harmonic analysis.

The final Section  \ref{S77} discusses the situation when the right-shift operator $S$
fails to map $X$ into itself.


\section{Preliminaries for sequence spaces}
\label{S2}


Recall that $\cno$ is a locally convex Fréchet space for the coordinate-wise convergence
topology. This topology is determined by the increasing sequence of seminorms
\begin{equation}\label{23}
r_n(x):=\max_{0\le k\le n}|x_k|,\quad x=(x_n)_\subi\in\cno,
\end{equation}
for each $n\in\no$. In $\cno$ we consider the 
coordinate-wise order for its positive cone. Namely,  
given $x=(x_n)_\subi\in\cno$ define $|x|:=(|x_n|)_\subi\in\cno$, then define $x\ge0$ if $x=|x|$ and finally define $x\le y$ if $(y-x)\ge0$. 
Accordingly, $\cno=\R^{\no}\oplus i\R^{\no}$ is the 
complexification of $\R^{\no}$.

A \textit{sequence space}  is a linear subspace  $X\subseteq\cno$. It is
a \textit{Banach sequence space} if it is endowed with a norm
for which it is complete.
A  \textit{Riesz norm} in $X$ is a norm $\|\cdot\|_X$ satisfying 
$\|x\|_X\le\|y\|_X$ whenever  $x,y\in X$ satisfy
$|x|\le|y|$. A  \textit{Banach lattice} $X\subseteq\cno$ is a 
Banach  sequence space  endowed with a Riesz norm for which it is a 
(complex) \textit{lattice}
for the coordinate-wise order. Namely,  for all $x,y\in X\cap \R^{\no}$ 
the supremum $x\vee y$ and
the infimum $x\wedge y$ exists in $ \R^{\no}$ and belong to $X$;
see, for example,  \cite[\S2.2]{MN} and \cite[Ch.II, \S11]{Sch}. 
A  Banach lattice $X\subseteq\cno$  is  \textit{solid} for the
coordinate-wise order if, whenever  $y\in \C^{\N_0}$ satisfies $|y|\le|x|$ with $x\in X$,
then $y\in X$. 
In this paper we deal with a particular class of sequence spaces.
Let $\ee:=\{e_n:n\in\no\}$, where $e_n:=(\delta_{nk})_{k=0}^\infty$ for each $n\in\no$;
it is an unconditional basis for $\cno$. 
We will say that $X\subseteq\cno$ is a \textit{solid Banach lattice},
briefly s.B.l.,  if it is a Banach lattice which is solid, 
contains $c_{00}:=\text{span}(\ee)$, and 
the natural inclusion map $X\subseteq\cno$ is continuous.
If, in addition, $\ee$ is a basis for $X$, then we say that $X$ is a 
s.B.l.\ with \textit{natural basis}. A s.B.l.\  
$X$ is called \textit{translation invariant} if the 
\textit{right-shift operator} $S\colon\cno\to\cno$, defined by
\begin{equation}\label{24}
Sx:=(0,x_0,x_1,x_2,\dots),\quad x\in\cno,
\end{equation}
satisfies $S(X)\subseteq X$.  
It turns out that each  space $X$ in the list
$$
c_0;\, \ell^p\text{ for } p\in[1,\infty];\, ces_q \text{ for } q\in\{0\}\cup (1,\infty]; \;
d_p\text{ for } p\in[1,\infty); N^p (1<p<\infty)
$$
is a s.B.l.; for $ces_p$, $d_p$ and $N^p$ see Sections \ref{S6} and \ref{S9}.
The  Banach lattice $c\subseteq\cno$ is \textit{not} 
a s.B.l.\  as it fails to be solid. 
The Banach  sequence space   $bv\subseteq \cno$ is not even a   lattice.
Indeed, $x=\big((-1)^n/(n+1)^2\big)_\subi$ and 
$|x|=\big(1/(n+1)^2\big)_\subi$ both belong to
$bv$ but $\|x\|_{\text{bv}}=(\pi^2/3)>2=\|\,|x|\,\|_{\text{bv}}$.
In Section  \ref{S77} we present an example of a s.B.l.\ $X$
with natural basis which fails to be translation invariant.


A useful tool for working with certain sequence spaces  is to view them 
within  the theory of Banach function spaces. 
To this aim we interpret elements of $\cno$  as functions $x\colon\no\to\C$
and consider  on $\no$ the counting measure 
$\mu\colon2^{\no}\to[0,\infty]$ specified by $\mu(\{n\})=1$ for each $n\in\no$. 
This  positive, $\sigma$-finite measure space is \textit{resonant}
in the sense of  \cite[Theorem II.2.7(ii)]{B-S}. 
Note that $L^0(\mu)=\cno$
since elements of $L^0(\mu)$ are finite $\mu$-a.e.\ and $\mu$
has no non-trivial null sets. 
A \textit{function norm} is a mapping $\rho\colon\mathcal{M}^+\to[0,\infty]$, where
$\mathcal{M}^+:=\{x\in\cno:x\ge0\}$, satisfying the following properties;
see  \cite[Definition I.1.1]{B-S}.
\begin{itemize}
\item[(P1)] $\rho(f)=0\iff f=0$ a.e.; $\rho(af)=a\rho(f),\; a\ge0$;
$\rho(f+g)\le\rho(f)+\rho(g)$.
\item[(P2)] $0\le g\le f$ a.e. $\Rightarrow \rho(g)\le \rho(f)$.
\item[(P3)] $0\le f_n\uparrow f$ a.e.$\Rightarrow \rho(f_n)\uparrow \rho(f)$.
\item[(P4)] $\mu(E)<\infty \Rightarrow \rho(\chi_E)<\infty$.
\item[(P5)] $\mu(E)<\infty \Rightarrow \int_Ef\,d\mu \le C_E \rho(f)<\infty$, for
some constant $0<C_E<\infty$.
\end{itemize}
The collection $X=X(\rho)$ of all  $x\in\cno$
such that $\rho(|x|)<\infty$ is called a \textit{Banach function space}
(briefly, B.f.s.) on $\no$. Define
\begin{equation}\label{41}
\|x\|_X:=\rho(|x|),\quad x\in X.
\end{equation}
The sequence space $(X,\|\cdot\|_X)$ is a Banach space 
which contains $c_{00}$ and 
the natural inclusion $X\subseteq\cno$ is continuous;
see   \cite[Theorems I.1.4 and I.1.6]{B-S}.
Moreover, $X$ is a 
Banach lattice relative to the Riesz norm $\|\cdot\|_X$ 
and  $X$ is \textit{solid} in $\cno$. Hence,
every B.f.s.\ is a s.B.l. The difference between
these two concepts lies in the Fatou property (P3), assumed in \cite{B-S} as part
of the definition of a function norm.
Indeed, given a s.B.l.\ $X\subseteq\cno$ it becomes a
B.f.s.\ on $\no$ (except for property (P3)) relative to the function norm
$\rho_X\colon\mathcal{M}^+\to[0,\infty]$ given by
$$
\rho_X(x):=
\begin{cases}
\| x\|_X,\quad x\in X\cap \mathcal{M}^+;
\\
\infty, \quad\text{otherwise}.
\end{cases}
$$


Concepts from the theory of B.f.s.' are relevant.
The distribution function $\mu_x\colon[0,\infty)\to[0,\infty]$ of $x\in\cno$
is given by
\begin{equation}\label{42}
\mu_x(\lambda):=\mu(\{n\in\no:|x_n|>\lambda\}),\quad \lambda\in[0,\infty);
\end{equation}
see   \cite[Definition II.1.1]{B-S}. Observe that $\mu_x$ assumes its values
in $\no\cup\{\infty\}$. Two elements $x,y\in\cno$ are called
\textit{equimeasurable} if $\mu_x=\mu_y$ (cf.   \cite[Definition II.1.2]{B-S}).
A B.f.s.\ 
$(X,\|\cdot\|_X)$  on $\no$,
with $\|\cdot\|_X$ given by  \eqref{41}, is \textit{rearrangement invariant} (briefly, r.i.)
if $\|x\|_X=\|y\|_X$ whenever
$x,y\in X$ are equimeasurable; see
 \cite[\S II.4]{B-S}. 
Every r.i.\ space $X$ over 
$\no$ satisfies
\begin{equation}\label{45}
\ell^1\subseteq X\subseteq\ell^\infty
\end{equation}
with both of the natural inclusion maps in \eqref{45}
having operator norm 1; see  \cite[Corollary II.6.8]{B-S}. 
Throughout this paper a B.f.s.\ $X=X(\rho)$ over $\no$ always
means relative to counting measure $\mu$ with the norm 
given by \eqref{41}.


Let $X$ be a B.f.s.\ on $\no$. 
An element $x\in X$ is said to have
\textit{absolutely continuous norm} (briefly, a.c.-norm) if
$\lim_{n\to\infty}\|x\chi_{E_n}\|_X=0$ for every sequence of sets
$\{E_n\}_{n=1}^\infty\subseteq 2^{\no}$ satisfying $E_n\downarrow\emptyset$; 
see   \cite[Definition I.3.1 and Proposition I.3.2]{B-S}. 
If every $x\in X$ has a.c.-norm, then $X$ is said to have a.c.-norm. 
The concept of absolute continuity is very important in the study of B.f.s.' The 
following  result  relates it to other useful concepts in our setting.

\begin{lemma}\label{l-41}
For a  B.f.s.\ $X$ on $\no$ the following conditions are equivalent.
\begin{itemize}
\item[(i)] $X$ has a.c.-norm.
\item[(ii)] $c_{00}$ is a dense subspace of  $X$.
\item[(iii)] $X$ is separable.
\item[(iv)] $\ee$ is a basis for $X$.
\end{itemize}
\end{lemma}

\begin{proof}
(i)$\iff$(iii). Since $\mu$ is a separable measure, apply
Corollary I.5.6 in \cite{B-S}.

(i)$\Rightarrow$(iv). Let $x=(x_n)_\subi\in X$ and define 
$x^{(N)}:=(x_0,x_1,\dots,x_N,0,0,\dots)=\sum_{k=0}^N x_ke_k$ for each $N\in\no$.
Clearly $|x^{(N)}|\le|x|$, for each $N\in\no$, and $\lim_{N\to\infty}x^{(N)}=x$
pointwise on $\no$. Accordingly, Proposition I.3.6 in \cite{B-S} 
implies that $\|x-x^{(N)}\|_X\to0$ for $N\to\infty$, that is, 
$x=\lim_{N\to\infty}\sum_{k=0}^N x_ke_k$ in $(X,\|\cdot\|_X)$. This shows
that the series $\sum_\subi x_ne_n$ converges to $x$ in $X$. Since $x\in X$ is arbitrary,
we can conclude that $\ee$ is a basis for $X$.

(iv)$\Rightarrow$(ii). Clear.

(ii)$\Rightarrow$(iii). The subset $Y\subseteq c_{00}$ consisting of all elements 
$x=(x_n)_\subi$ such that each coordinate $x_n$, for $n\in\no$, of $x$ is a rational 
(complex) number is countable. Clearly each element of $c_{00}$ is approximable in
$(X,\|\cdot\|_X)$ by a sequence of elements from $Y$. Since the closure
$\overline{c_{00}}=X$, it follows that $Y$ is a countable dense subset of $X$,
that is, $X$ is separable.
\end{proof}


We end this section by recalling some 
basic concepts from spectral theory. 
For a continuous linear operator $T\colon X\to X$, with $X$ a Banach space, 
its \textit{resolvent set} is defined by
$$
\rho(T;X):=\{\lambda\in\C: T-\lambda I \text{ is invertible in } \mathcal{L}(X)\},
$$
where $ \mathcal{L}(X)$ is the space of all continuous linear operators 
from $X$ into itself and $I\in  \mathcal{L}(X)$ is the identity operator. 
The operator norm of $T\in \mathcal{L}(X)$ is denoted by 
$\|T\|_{X\to X}$ and the dual Banach space of $X$ is denoted by $X'$.
The complement $\C\setminus\rho(T;X)$, denoted by $\sigma(T;X)$,
is called the \textit{spectrum} of $T$. The set $\sigma(T;X)$ is 
decomposed into three mutually disjoint parts which constitute the \textit{fine 
spectrum} of $T$. Namely, with $\mathcal{R}(S):=S(X)=\{Sx:x\in X\}$ denoting the
range of any $S\in \mathcal{L}(X)$, we have the \textit{point spectrum}
$$
\sigma_{\text{p}}(T;X):=\{\lambda\in\C: T-\lambda I \text{ is not injective}\}
$$
of $T$, the \textit{continuous spectrum}
$$
\sigma_{\text{c}}(T;X):=\{\lambda\in\C: T-\lambda I \text{ is injective
and $\mathcal{R}(T-\lambda I)$ is a proper, dense subspace of $X$}\}
$$
of $T$, and the \textit{residual spectrum}
$$
\sigma_{\text{r}}(T;X):=\{\lambda\in\C: T-\lambda I \text{ is injective
and the closure  $\overline{\mathcal{R}(T-\lambda I)}\not=X$}\}
$$ of $T$. Points in $\sigma_{\text{p}}(T;X)$ are called \textit{eigenvalues} of $T$.
The non-negative number  $r(T;X):=\sup\{|\lambda|:\lambda\in\sigma(T;X)\}$
is  the \textit{spectral radius} of $T$ and satisfies
$r(T;X)\le \|T\|_{X\to X}$. A standard reference for spectral theory is \cite{Co}, for example.


\section{$\ct$ and related operators: The main result}
\label{S3}


It is clear from \eqref{11}
that each operator $\ct$ is a linear map which is represented by a lower triangular matrix
with respect to the canonical, unconditional basis $\mathcal{E}$ of $\cno$. Namely,
\begin{equation}\label{21}
\ct \sim 
\begin{pmatrix}
1 & 0 & 0&0& \cdots \\
t/2 & 1/2 & 0&0& \cdots \\
t^2/3 & t/3 & 1/3&0& \cdots \\
t^3/4 & t^2/4 & t/4&1/4& \cdots \\
\cdot & \cdot & \cdot&\cdot& \cdot \\
\cdot & \cdot & \cdot&\cdot& \cdot 
\end{pmatrix}
,\quad t\in[0,1],
\end{equation}
with main diagonal $(1/(n+1))_{n=0}^\infty$.


\begin{lemma}\label{l-21}
Let $t\in[0,1]$.
\begin{itemize}
\item[(i)] Each $\ct$ is a positive operator on $\cno$, that is, $\ct x\ge0$ whenever $x\ge0$.
\item[(ii)] Let $0\le r\le s\le 1$. Then
\begin{equation}\label{22}
0\le |\mathcal{C}_r x|\le \mathcal{C}_r |x| \le \mathcal{C}_s |x|,
\quad x\in\cno.
\end{equation}

\item[(iii)] For each $t\in[0,1)$ we have the identities
$$
\ct e_n=\sum_{k=0}^\infty\frac{t^k}{k+n+1} e_{k+n}\in\ell^1,\quad n\in\no,
$$
and 
$$
\ct (e_n-t  e_{n+1})=\frac{1}{n+1} e_n,\quad n\in\no.
$$
In particular, the range of  $\ct$ contains $\ee$
whenever $\ct$ is defined on a sequence space  which contains $c_{00}$.

\item[(iv)] The linear map $\ct\colon\cno\to\cno$ is injective.
\end{itemize}
\end{lemma}


\begin{proof} 
(i) Immediate from the definition of $\ct$; see \eqref{11} and  \eqref{21}.

(ii) Fix $x\in\cno$. Let $0< r\le s\le 1$. For the $n$-th coordinate, 
\eqref{11} implies
$$
0\le \big(|\mathcal{C}_rx|\big)_n\le \frac{1}{n+1}\sum_{k=0}^nr^{n-k}|x_k|
=\big(\mathcal{C}_r|x|\big)_n
\le \frac{1}{n+1}\sum_{k=0}^ns^{n-k}|x_k|
=\big(\mathcal{C}_s|x|\big)_n,
$$
for each $n\in\no$, which is \eqref{22}. 
For $r=0$ and $s\in(0,1]$ it follows from \eqref{11} and \eqref{12} that
$$
0\le \big(|\mathcal{C}_0x|\big)_n=\frac{|x_n|}{n+1}
\le \frac{1}{n+1}\sum_{k=0}^ns^{n-k}|x_k|
=\big(\mathcal{C}_s|x|\big)_n,
$$
for each $n\in\no$, which is again \eqref{22}.

(iii) Both of the stated formulae follow by direct calculation from
the definition of $\ct$ (see also \eqref{21}) for $t\in[0,1)$ and 
the definition  of  $e_n$ for $n\in\no$.

(iv) Let $x\in\cno$ satisfy $\ct x=0$. It follows from \eqref{11} that
$$
t^nx_0+t^{n-1}x_1+\cdots+x_n=0,\quad n\in\no.
$$
The choice $n=0$ yields $x_0=0$. Using $x_0=0$ the choice $n=1$ yields
$x_1=0$. Using $x_0=x_1=0$ the choice $n=2$ yields
$x_2=0$. Proceeding inductively gives $x_n=0$ for all $n\in\no$, that is, $x=0$.
\end{proof}


Consider  the right-shift operator $S\colon\cno\to\cno$ given in \eqref{24}.
For each $n\in\N$ observe that 
$r_n(Sx)=\max_{0\le k< n}|x_k|\le r_n(x)$ 
and for $n=0$  that $r_0(Sx)=0\le r_0(x)$ for each $ x\in\cno$.
So, for every $n\in\no$, the operator $S$ satisfies
$r_n(Sx)\le r_n(x),\quad x\in\cno$,
which implies that $S\in\mathcal{L}(\cno)$. 
Note that each $r_n$, for $n\in\no$, is a \textit{Riesz seminorm} in $\cno$.
Accordingly, $\cno$ is a (complex) Fréchet lattice.


\begin{lemma}\label{l-22}
For each $t\in[0,1]$ the operator $\ct\colon\cno\to\cno$ is continuous. Moreover,
with $D$ being the diagonal operator in \eqref{12}, it is the case that
\begin{equation*}
\ct=\sum_\subi t^nDS^n=D\sum_\subi t^nS^n, \quad t\in[0,1),
\end{equation*}
where the series  converges for the strong operator topology 
in $\mathcal{L}(\cno)$, that is, 
\begin{equation}\label{25}
\ct x=\sum_\subi t^nDS^nx
\end{equation}
with the series being convergent  in $\cno$ for each $x\in\cno$.
\end{lemma}


\begin{proof}
Fix $t\in[0,1]$. Given $n\in\no$ it follows from \eqref{11} and \eqref{23} that
\begin{align*}
r_n(\ct x)&=\max_{0\le k\le n} \Big|\frac{1}{k+1}\sum_{j=0}^k t^{k-j}x_j\Big|
\le \max_{0\le k\le n} \frac{1}{k+1}\sum_{j=0}^k t^{k-j}|x_j|
\\ & \le \max_{0\le k\le n} \frac{1}{k+1}\sum_{j=0}^k |x_j|
\le r_n(x), 
\end{align*}
for each $x\in\cno$. Since $n\in\no$ is arbitrary, it follows that 
$\ct\in \mathcal{L}(\cno)$.

It is clear from \eqref{12} and the identity
\begin{equation}\label{26}
S^mx=(0,0,\dots0,x_0,x_1,x_2,\dots),\quad x\in \cno,\; m\in\no,
\end{equation}
with 0 occurring $m$-times before $x_0$, that 
\begin{equation}\label{27}
t^nDS^nx=\Big(0,0,\dots,0,\frac{t^nx_0}{n+1},\frac{t^nx_1}{n+1},
\frac{t^nx_2}{n+1},\dots\Big),\quad x\in\cno,
\end{equation}
for all $n\in\N$, with 0 occurring $n$-times prior to  $\frac{t^nx_0}{n+1}$. It is then
clear, algebraically, that \eqref{25} is valid by comparing the coordinates
in each side of \eqref{25}. Since $\ee$ is a  basis for $\cno$, it
remains to verify that the partial sums of the series in 
 \eqref{25} actually converge to the element $\ct x$ in the topology of $\cno$, for each $x\in\cno$.

So, fix $x\in\cno$. Select a fixed seminorm $r_k$ for some $k\in\no$. It is
clear from \eqref{11} and \eqref{27} that the coordinates in positions $0,1,\dots, k$ of
$\ct x- \sum_{j=0}^nt^jDS^jx$  are equal to $0$ for every $n>k$. 
Consequently, it follows that 
$r_k\big(\ct x- \sum_{j=0}^nt^jDS^jx\big)=0$, for all $n>k$,
which implies that $\lim_{n\to\infty}r_k\big(\ct x- \sum_{j=0}^nt^jDS^jx\big)=0$. Since
$k\in\no$ is arbitrary, we can conclude that the series in \eqref{25}
converges to $\ct x$ in $\cno$.
\end{proof}


\begin{lemma}\label{l-31}
Let $X$ be a s.B.l.\ and $D$ be the diagonal operator in 
\eqref{12}. Then $D\colon X\to X$  satisfies $\|D\|_{X\to X}=1$ and is 
a compact operator.  
\end{lemma}

\begin{proof}
Fix $x=(x_n)_\subi\in X$. Then $|Dx|=(|x_n|/(n+1))_\subi\le |x|$ and so $Dx\in X$ as $X$ is solid in $\cno$. Moreover,
$\|Dx\|_X=\|\,|Dx|\,\|_X\le \|\,|x|\,\|_X=\|x\|_X$.
Accordingly, $D\in\mathcal{L}(X)$ and $\|D\|_\txop\le1$.
Since $De_0=e_0$, it follows that
$\|De_0\|_X=\|e_0\|_X$ which, together with $\|D\|_\txop\le1$, implies that $\|D\|_\txop=1$.

Fix $n\in\no$ and define the finite rank operator $D^{[n]}\colon X\to X$ by 
$$
D^{[n]}x:=\Big(x_0,\frac{x_1}{2},\dots,\frac{x_n}{n+1},0,0,\dots\Big)=
\sum_{k=0}^n\frac{x_k}{k+1}e_k,\quad x\in X.
$$
Direct calculation yields that
$$
\big|(D-D^{[n]})x\big|=\Big(0,0,\dots,\frac{|x_{n+1}|}{n+2},\frac{|x_{n+2}|}{n+3},\dots\Big)
\le 
\frac{1}{n+2}|x|
$$
and hence,  that $\big\|(D-D^{[n]})x\big\|_X\le \frac{1}{n+2}\|x\|_X$, for $x\in X$.
This inequality implies that $\|D-D^{[n]}\|_{\txop}\le \frac{1}{n+2}$, for each $n\in\no$.
Accordingly, $\lim_{n\to\infty}D^{[n]}=D$ for the operator norm
in $\mathcal{L}(X)$, which implies  that $D$ is compact.
\end{proof}


Suppose now that  $X$ is  a translation invariant  s.B.l., 
that is, $S(X)\subseteq X$. Then
necessarily $S\in\mathcal{L}(X)$. Indeed, let $\{x^{[N]}\}_{N=1}^\infty\subseteq X$ 
be a sequence such that $\lim_{N\to\infty} x^{[N]}=0$ in $X$ and 
 $\lim_{N\to\infty}Sx^{[N]}=y$ in $X$ for some $y\in X$.
Since the natural inclusion map $X\subseteq \cno$ is continuous, it follows
that $x^{[N]}\to0$ in $\cno$ and $Sx^{[N]}\to y$ in $\cno$ for $N\to\infty$.
But, $S\in\mathcal{L}(\cno)$ and so $y=0$. Hence, $S\colon X\to X$ is a closed operator.
By  the Closed Graph Theorem we can conclude that $S\in\mathcal{L}(X)$.


\begin{proposition}\label{p-32}
Let $t\in[0,1)$ and $X$ be a translation invariant s.B.l.\  with natural basis. 
Suppose that the right-shift operator
$S\in\mathcal{L}(X)$ has the property that each series
\begin{equation}\label{32}
R_tx:=\sum_\subi t^nS^nx,\quad x\in X,
\end{equation}
is convergent in $X$. Then $\ct\in\mathcal{L}(X)$ and $\ct$ is a compact operator.
\end{proposition}

\begin{proof}
Since the sequence
$\big\{\sum_{n=0}^N t^nS^n\big\}_{N=0}^\infty\subseteq \mathcal{L}(X)$,
it follows from the Banach-Steinhaus Theorem that the pointwise defined linear operator
$R_t\colon x\mapsto R_tx$, for $x\in X$, is continuous, that is, $R_t\in\mathcal{L}(X)$.
Combined with Lemma \ref{l-31} we can conclude that 
$DR_t\in\mathcal{L}(X)$ is a compact operator. Moreover, since 
$D\in\mathcal{L}(X)$, the series expansion
\begin{equation}\label{33}
DR_tx=\sum_\subi t^nDS^nx
\end{equation}
is valid in $X$, for each $x\in X$. But, $X\subseteq\cno$ continuously, and so  we can 
conclude that the series in \eqref{33} also converges to $DR_tx$
in $\cno$, for each $x\in X$. Lemma \ref{l-22} shows that the series expansion in
\eqref{25}  is also valid in $\cno$. In particular, we can conclude that
$\ct x=DR_t x$, for $x\in X$.
But, $DR_t\in\mathcal{L}(X)$ and hence, also $\ct\in\mathcal{L}(X)$.
Since $DR_t\in\mathcal{L}(X)$ is compact, it follows that also $\ct$
is compact.
\end{proof}


We proceed to present a result which gives a positive solution
to the Conjecture in  the Introduction 
for a large class of sequence spaces $X$. 
The Conjecture requires the somewhat technical condition that the set
$\mathcal{S}$ (cf.  \eqref{31}) satisfies
$\mathcal{S}\subseteq X$. We will work with a more tractable substitute for this
condition.
The usual norm in the s.B.l.\  $\ell^p=\ell^p(\no)$, 
for $p\in[1,\infty]$, is denoted by $\|\cdot\|_p$.


Define the space
\begin{equation}\label{34}
d_1:=\Big\{x\in\ell^\infty:\hat x:=\Big(\sup_{k\ge n}|x_k|\Big)_\subi\in\ell^1\Big\};
\end{equation}
see \cite{Be}, \cite{curbera-ricker4}. 
The sequence $\hat x$ is called the \textit{least decreasing majorant} of $x$.
Then $d_1$ is a Banach lattice for the Riesz norm $\|x\|_{d_1}:=\|\hat x\|_1$ and the
coordinate-wise order. Since
\begin{equation}\label{35}
0\le |x|\le \hat x,\quad x\in\ell^\infty,
\end{equation}
it is clear that $\|x\|_1\le\|x\|_{d_1}$ for $x\in d_1$, that is, the natural inclusion
map $d_1\subseteq \ell^1$ has operator norm at most 1. But, 
$\|e_0\|_1=\|e_0\|_{d_1}$ and so the operator norm is precisely 1. 
The containment $d_1\subseteq \ell^1$ is \textit{proper}; see Remark \ref{r-610}(i)
below.

\begin{lemma}\label{l-33}
The set $\mathcal{S}\subseteq d_1$.
\end{lemma}

\begin{proof}
Fix $x\in\mathcal{S}$. Then
$|x|$ is \textit{eventually decreasing}  (as $\beta(x)<1$; see \eqref{31}), that is, 
there exists $N\in\no$ such that $(\hat x)_n=|x_n|$ for
all $n\ge N$; see the definition of $\hat x$ in \eqref{34}.
The ratio test implies, via \eqref{31}, that $x\in\ell^1$ and hence,  also 
$\hat x\in\ell^1$. Accordingly, $x\in d_1$. 
\end{proof}


\begin{lemma}\label{l-34}
Let $t\in[0,1)$. For each $m\in\no$ set $\lambda_m:=\frac{1}{m+1}$.
Define $x^{[0]}:=(t^n)_\subi$ and, for each $m\in\N$, define  $x^{[m]}\in\cno$  
coordinate-wise by 
$(x^{[m]})_j:=0$, for $0\le j<m$, with 
$(x^{[m]})_m=1$  and
\begin{equation}\label{36}
(x^{[m]})_{m+n}:=\frac{(m+1)(m+2)\cdots(m+n)}{n!}\,t^n,\quad n\ge1.
\end{equation}
\begin{itemize}
\item[(i)] Let $m\in\no$. Then $\lambda_m$ is an eigenvalue of
$\ct\in\mathcal{L}(\cno)$ and its corresponding eigenspace
\begin{equation}\label{36-37}
\Big\{x\in\cno: \ct x=\lambda_m x\Big\}=
\Big\{\alpha x^{[m]}:\alpha\in\C\Big\}
\end{equation}
is one-dimensional.
\item[(ii)] The subset $\big\{x^{[m]}:m\in\no\big\}$ 
of $\cno$ is actually contained in $d_1$.
\end{itemize}
\end{lemma}

\begin{proof}
(i) For the details of the purely algebraic calculations required to establish  \eqref{36-37} we refer to the proof of Theorem 9 in \cite{Y-M-D}, for example.
These calculations are 
carried out in $\cno $ and only afterwards  is it shown that $x^{[m]}$ lies in
the relevant space being considered in \cite{Y-M-D}.

(ii) Consider $m=0$. Then $|x^{[0]}|=(t^n)_\subi$ is a decreasing 
sequence and so 
$(x^{[0]})\hat{}=|x^{[0]}|\in\ell^1$. Accordingly, $x^{[0]}\in d_1$.

Now fix $m\ge1$ and $t\in(0,1)$. 
According to \eqref{36} the ratio
\begin{equation}\label{37}
\frac{(x^{[m]})_{m+n+1}}{(x^{[m]})_{m+n}}=\frac{(m+n+1)t}{n+1},\quad n\ge1,
\end{equation}
is at most 1 whenever $n\ge (t(m+1)-1)/(1-t)$. So, there exists a smallest
integer $N(t,m)\in\no$ such that
$0\le (x^{[m]})_{m+n+1} \le (x^{[m]})_{m+n}$, for $n\ge N(t,m)$,
that is, the sequence $x^{[m]}$ is eventually decreasing. Hence,
$x^{[m]}\in d_1$ if and only if $x^{[m]}\in \ell^1$. But, 
\eqref{37} implies that
$$
\lim_{n\to\infty} \frac{(x^{[m]})_{m+n+1}}{(x^{[m]})_{m+n}}=t\in(0,1)
$$
from which it follows that $0\le x^{[m]}\in \ell^1$, that is, $x^{[m]}\in d_1$.

For the case $t=0$ we see that $x^{[m]}= e_m\in d_1$.
\end{proof}


The main result of the paper   can now be presented.

\begin{theorem}\label{t-35}
Let $t\in[0,1)$ and $X$ be a translation invariant s.B.l.\
with  natural basis which contains $d_1$. 
Suppose that  the right-shift operator $S\in\mathcal{L}(X)$ has the property
 that the  series in \eqref{32} is convergent in $X$ for each $x\in X$.
Then  $\ct\in\mathcal{L}(X)$  is a compact operator with spectrum
\begin{equation}\label{38}
\sigma(\ct;X)=\Lambda\cup\{0\}. 
\end{equation}
Moreover, the fine spectra of $\ct$ consist of the point spectrum
\begin{equation}\label{39}
\sigma_{\text{p}}(\ct;X)=\Lambda 
\end{equation}
together with
\begin{equation}\label{310}
\sigma_{\text{c}}(\ct;X)=\{0\}\quad\text{and}\quad \sigma_{\text{r}}(\ct;X)=\emptyset.
\end{equation}
\end{theorem}

\begin{proof}
Since  $X$ has natural basis, the operator $\ct\in\mathcal{L}(X)$
is lower triangular with entries $1/(n+1)$, for $n\in\no$ on the main diagonal; see
\eqref{21}. By applying Lemma 2.2 of \cite{S-El} we can conclude that
\begin{equation}\label{311}
\sigma_{\text{p}}(\ct;X)\subseteq\Lambda. 
\end{equation}
Lemma \ref{l-34} shows, for each $m\in\no$, that there exists $x\in d_1\setminus\{0\}$
satisfying $\ct x=\frac{1}{m+1}x$. Since $d_1\subseteq X$ (by assumption), 
it follows that $\frac{1}{m+1}\in \sigma_{\text{p}}(\ct;X)$ and so 
 \eqref{311} is an equality, that is, \eqref{39} is valid.

In view of \eqref{39} and the fact that $\ct\in\mathcal{L}(X)$ is compact,
the classical result of F. Riesz for compact operators implies that
\eqref{38} is valid; see \cite[Theorem VII.7.1]{Co}.

Since $\ct\in\mathcal{L}(X)$ is injective (cf.\ Lemma \ref{l-21}(iv)), 
0 belongs to either the residual
or the continuous spectrum of $\ct$. By Lemma \ref{l-21}(iii)  the
canonical vectors satisfy $\ee\subseteq\ct(X)\subseteq X$.
Since $X$ has natural basis, it follows that the range
$\ct(X)$ is dense in $X$. Hence, $0\in \sigma_{\text{c}}(\ct;X)$
from which \eqref{310} follows.
\end{proof}


\begin{remark}\label{r-36}
The natural inclusion map $j\colon d_1\to X$,
whenever $d_1\subseteq X$ with $X$ a s.B.l., is necessarily
continuous. Indeed, let  $\{x^{[N]}\}_{N=1}^\infty\subseteq d_1$ be a 
sequence
such that $\lim_{N\to\infty}x^{[N]}=0$ in $d_1$ and 
$\lim_{N\to\infty}j(x^{[N]})=y$ in $X$ for some $y\in X$. Since 
$d_1\subseteq \ell^1$ continuously (see the discussion prior to
Lemma \ref{l-33}) with $\ell^1\subseteq \cno$ continuously and  $X\subseteq\cno$ continuously (by assumption),
it follows that $x^{[N]}\to0$ in $\cno$ and $x^{[N]}=j(x^{[N]})\to y$ in $\cno$ for $N\to\infty$.
This implies that $y=0$. Hence,  $j$ is a closed operator. By the Closed Graph Theorem
$j$ is continuous.
\end{remark}

The remainder of the paper is devoted to applications of Theorem
\ref{t-35}.  We begin immediately with a large class of classical spaces.


Let $X$ be a r.i.\ space over $\no$. 
Let us review the  assumptions of Theorem \ref{t-35} for such spaces.
Since $d_1\subseteq \ell^1$, it follows from \eqref{45} that
$d_1\subseteq X$. 
It is clear from \eqref{42} that 
\begin{equation*}
\mu_x=\mu_{Sx},\quad x\in\cno,
\end{equation*}
which implies that 
$\|Sx\|_X=\|x\|_X$, for $x\in X$,
that is, $S\in\mathcal{L}(X)$ is an isometry. 
In particular,  $X$ is translation invariant and 
$S$ satisfies
\begin{equation}\label{44}
\sum_\subi t^n\|S^n\|_{X\to X}=\sum_\subi t^n<\infty,\quad t\in[0,1).
\end{equation}
Accordingly, the series for $R_t$ given in \eqref{32} is absolutely
convergent for the operator norm (hence, also for the strong operator topology) in
$\mathcal{L}(X)$.
So, $X$  satisfies all of the assumptions of Theorem \ref{t-35} except for one.
Namely, although the set $\ee\subseteq X$, it need \textit{not} necessarily
be a basis for $X$. Consider, for example,  $X=\ell^\infty$. 
Actually, $X$ has natural basis
if and only if it is separable; see Lemma \ref{l-41}.


The above discussion, for the setting of the following result, shows that all 
the assumptions of Theorem \ref{t-35} are fulfilled. Note, by Lemma \ref{l-31}
and \eqref{44}, together with $\ct=DR_t$ (see the proof of Proposition \ref{p-32}),
that 
$$
\|\ct\|_{\txop}\le\|D\|_{\txop}\|R_t\|_{\txop}
\le \sum_{n=0}^\infty t^n=\frac{1}{1-t}.
$$


\begin{theorem}\label{t-43} 
Let $X$ be a separable, r.i.\ space over $\no$
and $t\in[0,1)$.
The  operator $\ct\in\mathcal{L}(X)$ satisfies 
\begin{equation}\label{46}
\|\ct\|_{\txop}\le (1-t)^{-1},
\end{equation}
is compact  and has spectrum
$$
\sigma\big(\ct; X\big)=\Lambda\cup\{0\}.
$$
Moreover, the fine spectra of $\ct$ are given by
$$
\sigma_{\text{p}}\big(\mathcal{C}_t; X\big)=\Lambda,\quad
\sigma_{\text{c}}\big(\mathcal{C}_t; X\big)=\{0\},\quad
\sigma_{\text{r}}\big(\mathcal{C}_t; X\big)=\emptyset.
$$
\end{theorem}


\begin{remark}\label{r-44}
(i) Observe that the right-side of  \eqref{46} is independent of the r.i.\ space $X$.

(ii) The class of  r.i.\ spaces over $\no$  is quite extensive.
As a sample we mention that the spaces $\ell^{p,r}$ with $p,r\in[1,\infty)$ are separable
r.i.\ spaces; see Section IV.4 in \cite{B-S}. Note that $\ell^{p,p}=\ell^p$, for $p\in[1,\infty)$.
The Lorentz-Zygmund spaces $\ell^{p,r}(\log l)^\alpha$ are also r.i.\
(cf. pp.\ 284-285 in \cite{B-S}). The same is true for 
Orlicz sequence spaces $X\subseteq\cno$
which have the Fatou-property (cf. Chapter 4 in Volume I of \cite{LT} and
Chapter 2 in Volume II of \cite{LT}). If $X$ satisfies a 
$\Delta_2$-condition at 0 (cf. \cite[Vol.\ I, Definition 4.a.3]{LT}), 
then $X$ is separable; 
see \cite[Vol.\ I, Proposition 4.a.4]{LT}.
\end{remark}


\section{The operators $\ct$ acting in discrete Ces\`{a}ro  spaces and  their
duals}
\label{S6}


For $1<p<\infty$, the sequence spaces
\begin{equation}\label{51}
ces_p:=\Big\{x\in \C^{\N_0}:\|\mathcal{C}_1|x|\|_p<\infty\Big\},
\end{equation}
which are Banach lattices for the order from $\cno$ and when equipped with the 
Riesz norm
\begin{equation}\label{52}
\|x\|_{ces_p}:=\|\mathcal{C}_1|x|\|_p
=\left(\sum_{n=0}^\infty
\bigg(\frac{1}{n+1}\sum_{k=0}^n|x_k|\bigg)^p\right)^{1/p},\quad x\in ces_p\,,
\end{equation}
were thoroughly investigated by 
G. Bennett in the monograph \cite{Be}. According to Proposition 2 of \cite{Ja}
each space $ces_p$ is reflexive. In particular,
$ces_p$ has a.c.-norm. Moreover, $ces_p$ is  $p$-concave and 
$\ee$ is an unconditional basis, \cite[Proposition 2.1]{curbera-ricker4}. 
It is clear from \eqref{52} that $ces_p$ is solid.

The sequence space $ces_0$ is given by
\begin{equation}\label{53}
ces_0:=\Big\{x\in\cno:\mathcal{C}_1|x|\in c_0\Big\}.
\end{equation}
It is a  s.B.l.\ for the order from $\cno$ and for the Riesz norm 
\begin{equation*}
\|x\|_{ces_0}:=\|\mathcal{C}_1|x|\|_{c_0}=\|\mathcal{C}_1|x|\|_\infty,\quad x\in ces_0,
\end{equation*}
and this norm is a.c. Hence, the K\"othe dual of $ces_0$ 
coincides with its dual Banach space $ces_0'$ and 
we have that $ces_0'=d_1$. Moreover, $\ee$ is an unconditional basis
in $ces_0$. In particular, $ces_0$ is \textit{separable}. For these facts we
refer to Section 6 of \cite{curbera-ricker4}, for example.

The sequence space $ces_\infty$ is defined by
\begin{equation}\label{55}
ces_\infty:=\Big\{x\in\cno:\mathcal{C}_1|x|\in \ell^\infty\Big\}.
\end{equation}
It is a Banach lattice for the order from $\cno$ and for the Riesz norm
\begin{equation}\label{56}
\|x\|_{ces_\infty}:=\|\mathcal{C}_1|x|\|_\infty,\quad x\in ces_\infty.
\end{equation}
In Proposition \ref{p-76} it will be shown that $ces_\infty$ is \textit{non-separable}.
In particular, $\ee$ cannot be a basis for $ces_\infty$. 
The Banach space bidual $ces_0''=ces_\infty$; see 
\cite{Al} and \cite[Section 6]{curbera-ricker4}.


Given $p\in(1,\infty]$,
define a function norm $\rho_p\colon\mathcal{M}^+\to[0,\infty]$ via
$$
\rho_p(x):=\|\ce_1x\|_p\,,\quad x\in \mathcal{M}^+;
$$
see Section \ref{S2}. Recall that $\ce_1\colon\cno\to\cno$ is a positive operator
(cf.  Lemma \ref{l-21}(i)). In the notation of \eqref{41} we note that
$$
ces_p=X(\rho_p):=\{x\in\cno:\rho_p(|x|)<\infty\}
$$
is a B.f.s.\ on $\no$ with $\rho_p(x)=\|x\|_{ces_p}$
for $x\in ces_p$.
The space $ces_0=(ces_\infty)_a$ is the closed (proper) ideal in  $ces_\infty$
consisting of all the elements having a.c.-norm. 
Since the s.B.l.\ $ces_0$ fails to have the Fatou property,
it is not a B.f.s.\ on $\no$.


The Ces\`{a}ro  operator 
$ \mathcal{C}_1\colon ces_p\to ces_p$ is continuous with operator norm
$\| \mathcal{C}_1\|_{ces_p\to ces_p}=p'=p/(p-1)$ and its spectrum
$$
\sigma\big(\mathcal{C}_1; ces_p\big)=
\Big\{z\in\C:|z-\frac{p'}{2}|\le \frac{p'}{2}\Big\},\quad 1<p<\infty,
$$
where $\frac{1}{p}+\frac{1}{p'}=1$, \cite[Theorem 5.1]{curbera-ricker4}. 
Concerning the fine 
spectrum of $\mathcal{C}_1$ we note that 
$\sigma_{\text{p}}\big(\mathcal{C}_1; ces_p\big)=\emptyset$
and  $\sigma_{\text{r}}\big(\mathcal{C}_1; ces_p\big)
=\big\{z\in\C:|z-\frac{p'}{2}|< \frac{p'}{2}\big\}$;
see Theorem 5.1 of \cite{curbera-ricker4} and its proof. 
Consequently, 
$\sigma_{\text{c}}\big(\mathcal{C}_1; ces_p\big)
=\big\{z\in\C:|z-\frac{p'}{2}|= \frac{p'}{2}\big\}$.
So, the fine spectra of 
$\mathcal{C}_1\in\mathcal{L}(ces_p)$ are completely identified.
Since $\sigma\big(\mathcal{C}_1; ces_p\big)$ is an uncountable set, the operator
$\mathcal{C}_1\colon ces_p\to ces_p$ is  not compact.


It is known that $\mathcal{C}_1\colon ces_0\to ces_0$ is continuous with operator 
norm 1 and that
$$
\sigma\big(\mathcal{C}_1; ces_0\big)
=\Big\{z\in\C:|z-\frac{1}{2}|\le \frac{1}{2}\Big\},
$$
\cite[Theorem 6.4]{curbera-ricker4}. The proof of Theorem 
6.4 in  \cite{curbera-ricker4} reveals that
$\sigma_{\text{p}}\big(\mathcal{C}_1; ces_0\big)=\emptyset$, 
that $\sigma_{\text{r}}\big(\mathcal{C}_1; ces_0\big)=
\{1\}\cup\big\{z\in\C:|z-\frac{1}{2}|< \frac{1}{2}\big\}$ and that
$\sigma_{\text{c}}\big(\mathcal{C}_1; ces_0\big)=
\big\{z\in\C:|z-\frac{1}{2}|= \frac{1}{2}\big\}\setminus\{1\}$.

In contrast to  $\mathcal{C}_1$, the behavior of the generalized Ces\`{a}ro  
operators $\ct$ on
$ces_p$, for $t\in[0,1)$, is significantly different.
We first collect together a few basic facts.

\begin{lemma}\label{l-51}
Let $t\in[0,1]$ and $p\in(1,\infty]$.
\begin{itemize}
\item[(i)] The
operator $\ct\colon ces_p\to\ell^p$ satisfies $\|\ct\|_{ces_p\to \ell^p}\le1$. 
\item[(ii)] The
operator $\ct \colon ces_0 \to c_0$ satisfies $\|\ct\|_{ces_0\to c_0}=1$.
\item[(iii)] For $t\not=1$, the norm of $\ct\colon \ell^\infty\to\ell^p$ satisfies $\|\ct\|_{\ell^\infty\to\ell^p}\le\|\xi\|_p/(1-t)$, where
$\xi:=(1/(n+1))_{n=0}^\infty$.
\end{itemize}
\end{lemma}


\begin{proof}
(i) Fix $t\in[0,1]$ and $p\in(1,\infty]$. Set $r:=t$ and $s=1$ in  \eqref{22}  and recall that 
$\|\cdot\|_p$ is a  Riesz norm yields, for  $x\in ces_p$
(i.e., $\mathcal{C}_1|x|\in\ell^p$; see \eqref{51} and  \eqref{22}), 
that
$$
\|\ct x\|_p= \||\ct x|\|_p\le \|\ct |x|\|_p\le  
\|\mathcal{C}_1|x|\|_p=\|x\|_{ces_p},
$$
where the last equality is due to  \eqref{52} and \eqref{56}. It follows that
$\|\ct\|_{ces_p\to \ell^p}\le 1$.

(ii) Fix $t\in[0,1]$. For each $x\in ces_0$ (i.e., $\mathcal{C}_1|x|\in c_0$; see \eqref{53})
it follows from   \eqref{22} with $r:=t$ and $s=1$ that
$$
\|\ct x\|_{c_0}= \|\ct x\|_\infty\le \|\ct |x|\|_\infty\le  
\|\mathcal{C}_1|x|\|_\infty=\|x\|_{ces_0}.
$$
Since $x\in ces_0$ is arbitrary, it follows that 
$\|\ct\|_{ces_0\to c_0}\le1$. But, $\|e_0\|_{ces_0}=1$ and $\|\ct e_0\|_{c_0}=1$
and so actually $\|\ct\|_{ces_0\to c_0}=1$.

(iii) Fix $t\in[0,1)$ and $p\in(1,\infty]$. Let $x\in\ell^\infty$. For the $n$-th
coordinate we have
\begin{equation}\label{57}
\big(\ct|x|\big)_n=\frac{1}{n+1}\sum_{k=0}^n t^{n-k}|x_k|
\le
\frac{\|x\|_\infty}{n+1}\sum_{k=0}^n t^k 
\le\frac{\|x\|_\infty}{(n+1)(1-t)},\quad n\in\no.
\end{equation}
Since $\xi\in\ell^p$ for every $p\in(1,\infty]$, it follows that
$$
\|\ct x\|_p\le \|\ct|x|\|_p\le \frac{\|\xi\|_p}{(1-t)}\|x\|_\infty,\quad x\in\ell^\infty.
$$
In particular, $\ct\colon\ell^\infty\to\ell^p$ has operator norm at most
$\|\xi\|_p/(1-t)$.
\end{proof}


\begin{remark}\label{r-52} 
Concerning  Lemma \ref{l-51}(iii) when $p=\infty$ we note,
for $t=1$, that \eqref{57} becomes
$$
\big(\mathcal{C}_1|x|\big)_n\le \|x\|_\infty,\quad x\in\ell^\infty,\; n\in\no,
$$
because $\frac{\|x\|_\infty}{(n+1)}\sum_{k=0}^n1^k=\|x\|_\infty$. Accordingly,
$\mathcal{C}_1\colon\ell^\infty\to\ell^\infty$ with 
$\|\mathcal{C}_1\|_{\ell^\infty\to\ell^\infty}\le1$.
But, $\|e_0\|_{\infty}=1$ and $\|\mathcal{C}_1 e_0\|_{\infty}=1$
and so actually $\|\ct\|_{\ell^\infty\to\ell^\infty}=1$.
\end{remark}


\begin{lemma}\label{l-53} 
For  $p\in\{0\}\cup(1,\infty]$ the natural inclusion
$ces_p\subseteq \cno $ is continuous.
\end{lemma}

\begin{proof}
Fix $1<p<\infty$. According to Lemma 4.7 in \cite{Be} there exists $\alpha_p>0$ such that
$$
\frac{\alpha_p}{(j+1)^{1/p'}}\le \|e_j\|_{ces_p},\quad j\in\no,
$$ 
where $\frac1p+\frac{1}{p'}=1$. Fix $n\in\no$. Given $x\in ces_p$ it follows from
$|x_ke_k|\le |x|$  that
$$
\frac{\alpha_p|x_k|}{(k+1)^{1/p'}}\le \|x_ke_k\|_{ces_p}\le \|x\|_{ces_p}\,,
\quad 0\le k\le n,
$$ 
and hence, that
$$
|x_k|\le \alpha_p^{-1} (k+1)^{1/p'}\|x\|_{ces_p}
\le \alpha_p^{-1} (n+1)^{1/p'}\|x\|_{ces_p}.
$$
So, in the notation of \eqref{23} and with $A_n:=\alpha_p^{-1} (n+1)^{1/p'}$, we have
that $r_n(x)\le A_n \|x\|_{ces_p}$, for $x\in ces_p$.
Since $n\in\no$ is arbitrary, we can conclude that $ces_p\subseteq\cno$ continuously.

Let $p=0$. Then $\|e_n\|_{ces_0}=\|e_n\|_{ces_\infty}=\frac{1}{n+1}$ for all
$n\in\no$. Fix $n\in\no$. Given $x\in ces_0$ it follows, for each $0\le k\le n$, that
$$
|x_k|= (k+1)\|x_k e_k\|_{ces_0}
\le (k+1) \|x\|_{ces_0}\le (n+1) \|x\|_{ces_0}
$$
and hence, that $r_n(x)\le (n+1) \|x\|_{ces_0}$, for $x\in ces_0$.
The same calculations are valid in $ces_\infty$. Since $n\in\no$ is arbitrary, it
follows that both inclusions $ces_0\subseteq\cno$ and $ces_\infty\subseteq\cno$
are continuous.
\end{proof}


\begin{lemma}\label{l-54} 
Let $p\in\{0\}\cup(1,\infty]$. The right-shift operator
satisfies $\|S\|_{ces_p\to ces_p}\le1$.
\end{lemma}

\begin{proof}
Fix  $p\in(1,\infty)$. Given $x\in ces_p$ observe that 
$|Sx|=S|x|$ and hence, 
\begin{align*}
\|Sx\|_{ces_p}&= \|\mathcal{C}_1S|x|\|_p=
\Big\|\Big(0,\frac{|x_0|}{2}, \frac{|x_0|+|x_1|}{3},\frac{|x_0|+|x_1|+|x_2|}{4},\dots\Big)\Big\|_p
\\ & \le 
\Big\|\Big(0,|x_0|, \frac{|x_0|+|x_1|}{2},\frac{|x_0|+|x_1|+|x_2|}{3},\dots\Big)\Big\|_p
\\ & = 
\|\mathcal{C}_1|x|\|_p=\|x\|_{ces_p}.
\end{align*}
It follows that $\|S\|_{ces_p\to ces_p}\le1$.
Analogous calculations are valid for $p\in\{0,\infty\}$.
\end{proof}


\begin{remark}\label{r-55} 
It follows from Lemma \ref{l-31}, Lemma \ref{l-53} and \eqref{32} that
$$
\|\ct\|_{ces_p\to ces_p}=\|DR_t\|_{ces_p\to ces_p}\le
\|D\|_{ces_p\to ces_p}\|R_t\|_{ces_p\to ces_p}\le\sum_{n=0}^\infty t^n=\frac{1}{1-t},
$$
for every $p\in\{0\}\cup(1,\infty)$ and $t\in(0,1)$.
\end{remark}


Clearly each space $ces_p$, for $p\in\{0\}\cup(1,\infty]$, is a  s.B.l.\  
because of Lemma  \ref{l-53}
and the fact that $c_{00}\subseteq ces_p$. 
For $p\in(1,\infty]$ this also follows from the fact that $ces_p$ is 
a B.f.s.\ on $\no$; see Section \ref{S2}.
Moreover, Lemma \ref{l-54} shows that $ces_p$ is translation invariant and the right-shift
operator $S$ satisfies
$$
\sum_{n=0}^\infty \|t^nS^n\|_{ces_p\to ces_p}\le \sum_{n=0}^\infty t^n<\infty,
\quad t\in(0,1),\quad p\in\{0\}\cup(1,\infty].
$$
In particular, the series for  $R_t$ 
in \eqref{32} converges in $\mathcal{L}(ces_p)$ for the operator norm
(hence,  for the strong operator) topology. To verify  $d_1\subseteq ces_p$ fix
$p\in(1,\infty)$. Given $x\in d_1$, it follows from \eqref{35} that 
$\mathcal{C}_1|x|\le \mathcal{C}_1\hat x$. Since $\hat x\in\ell^1$ and 
$\mathcal{C}_1\in \mathcal{L}(\ell^p)$, we have  $\mathcal{C}_1\hat x\in \ell^p$
and hence, also $\mathcal{C}_1|x|\in \ell^p$. That is, $x\in ces_p$. Noting that 
$\hat x\in\ell^1\subseteq c_0\subseteq \ell^\infty$ and $\mathcal{C}_1$ is
continuous in both $c_0$ and $\ell^\infty$, the previous argument also shows that 
$\mathcal{C}_1|x|\in c_0\subseteq \ell^\infty$. That is,
$d_1\subseteq ces_0\subseteq ces_\infty$. With the exception of $ces_\infty$,
where $\ee$ is \textit{not} a basis, all the assumptions in Theorem \ref{t-35} 
have been verified. So, for $X=ces_p$,
we can deduce the following result (also Remark \ref{55}  is
used).


\begin{theorem}\label{t-56}
Let $t\in[0,1)$ and $p\in\{0\}\cup(1,\infty)$. 
The operator $\ct\in\mathcal{L}(ces_p)$
satisfies  
\begin{equation}\label{58}
\|\ct\|_{ces_p\to ces_p}\le (1-t)^{-1},
\end{equation}
is compact and  has spectrum
$$
\sigma\big(\ct; ces_p\big)=\Lambda\cup\{0\}.
$$
Moreover, the fine spectra of $\ct$ are given by
$$
\sigma_{\text{p}}\big(\mathcal{C}_t; ces_p\big)=\Lambda,\quad
\sigma_{\text{c}}\big(\mathcal{C}_t; ces_p\big)=\{0\},\quad
\sigma_{\text{r}}\big(\mathcal{C}_t; ces_p\big)=\emptyset.
$$
\end{theorem}


\begin{remark}\label{r-57} 
(i) An alternate estimate to \eqref{58} is possible. Fix $p\in(1,\infty)$ and $t\in[0,1)$.
According to Theorem 326 of \cite{HLP} the operator 
$\mathcal{C}_1\in \mathcal{L}(\ell^p)$ 
satisfies $\|\ce_1\|_{\ell^p\to\ell^p}=p/(p-1)=p'$.
Fix $x\in ces_p$. Lemma \ref{l-51}(i) implies that $\ct x\in\ell^p$
and hence, also $|\ct x|\in\ell^p$. So,
$$
\|\ct x\|_{ces_p}:=\|\mathcal{C}_1|\ct x|\|_p\le \frac{p}{p-1} 
\||\ct x|\|_p=\frac{p}{p-1} \|\ct x\|_p.
$$
But, $\|\ct x\|_p\le \|x\|_{ces_p}$ (cf. Lemma \ref{l-51}(i)) and hence,
$\|\ct x\|_{ces_p}\le \frac{p}{p-1}\|x\|_{ces_p}$.
Since $x\in ces_p$ is arbitrary, it follows that $\|\ct\|_{ces_p\to ces_p}\le \frac{p}{p-1}$.

Consider now $p=0$ and fix $x\in ces_0$. Then $\ct x\in c_0$ (cf. Lemma \ref{l-51}(ii))
and so, by Remark \ref{r-52}, we have that
$$
\|\ct x\|_{ces_0}:=\|\mathcal{C}_1|\ct x|\|_{c_0}\le \|\mathcal{C}_1\ct |x|\|_\infty
\le \|\mathcal{C}_1\|_{\ell^\infty \to \ell^\infty} \| \ct |x|\|_\infty=\| \ct |x|\|_\infty.
$$
But,  \eqref{22} implies that
$\| \ct |x|\|_\infty\le \| \mathcal{C}_1 |x|\|_\infty=\|x\|_{ces_0}$
from which it follows that $\|\ct x\|_{ces_0}\le \|x\|_{ces_0}$ for $x\in ces_0$
and so $\|\ct\|_{ces_0\to ces_0}\le1$. Since $\|e_0\|_{ces_0}=1$ 
and $\|\ct e_0\|_{ces_0}=\|\mathcal{C}_1\big(\frac{t^n}{n+1}\big)_\subi\|_\infty=1$
we actually have $\|\ct\|_{ces_0\to ces_0}=1$, for $t\in[0,1)$.
A similar calculation shows that
$\|\ct\|_{ces_\infty\to ces_\infty}=1$, for $t\in[0,1)$.

(ii) According to \eqref{45} the B.f.s.' $ces_p$ for $p\in(1,\infty)$
are \textit{not} r.i., due to the fact that they are not contained in $\ell^\infty$;
see Remark 2.2(ii) in \cite{curbera-ricker4}. It is routine
to check that $ces_p\subseteq ces_0$ for all $p\in(1,\infty)$ and
so the s.B.l.\  $ces_0$ cannot be contained in $\ell^\infty$ either.

(iii) It follows from \eqref{58} and part (i) of this Remark that
$$
\|\ct\|_{ces_p\to ces_p}\le\min\Big\{ \frac{p}{p-1},\frac{1}{1-t}\Big\},
\quad t\in[0,1),\; p\in(1,\infty).
$$
It is routine to verify that
$$
\|e_0\|_{ces_p}=\|\ce_1e_0\|_p<
\Big\|\ce_1\Big(\frac{t^n}{n+1}\Big)_\subi\Big\|_p=
\big\|\ce_1\ct e_0\big\|_p=\big\|\ct e_0\big\|_{ces_p},
\quad p\in(1,\infty),
$$
which implies that $\|\ct\|_{ces_p\to ces_p}>1$. According to  Theorem \ref{t-56}(ii)
the spectral radius $r(\ct;ces_p)=1$ and so we can conclude (unlike for $t=1$)
that
$$
r(\ct;ces_p)<\|\ct\|_{ces_p\to ces_p},\quad t\in(0,1),\; p\in(1,\infty).
$$
For $p\in\{0,\infty\}$ we have that
$r(\ct;ces_p)=\|\ct\|_{ces_p\to ces_p}=1$, for $t\in(0,1)$.
\end{remark}


Consider now the B.f.s.\  $ces_\infty$  given by
\eqref{55} and equipped with the Riesz norm
\eqref{56}. Recall that the Banach space bidual 
$ces_0''=ces_\infty$ and that the dual Banach space $ces_0'=d_1$. The duality is specified by $\langle u,x\rangle:=\sum_\subi u_nx_n$, for $u\in ces_0$ and  $x\in d_1$.

\begin{lemma}\label{l-75}
The natural inclusion $\ell^\infty\subseteq ces_\infty$ has operator norm equal to 1.
\end{lemma}

\begin{proof}
For $x\in\ell^\infty$ we have that
$|\mathcal{C}_1x|\le \mathcal{C}_1|x|\le \|x\|_\infty \mathbbm{1}$,
where $\mathbbm{1}:=(1,1,\cdots)$, and so
$\|x\|_{ces_\infty}=\|\mathcal{C}_1|x|\|_\infty\le \|x\|_\infty$
because $\|\mathcal{C}_1\mathbbm{1}\|_\infty=1$. Accordingly,
the operator norm of the inclusion is
at most 1. That it actually equals 1 follows from 
$\|\mathbbm{1}\|_\infty=\|\mathbbm{1}\|_{ces_\infty}=1$.
\end{proof}

\begin{proposition}\label{p-76}
\begin{itemize}
\item[(i)] The B.f.s.\  $ces_\infty$ is non-separable.
\item[(ii)]  $\ell^\infty$ is not dense in $ces_\infty$.
\end{itemize}
\end{proposition}

\begin{proof}
(i) Consider the sequence in $ces_\infty$ given 
by $\mathbbm{1}^{[N]}:=\sum_{k=0}^{N-1} e_k$
for $N\in\N$. Then $0\le \mathbbm{1}^{[N]} \le \mathbbm{1}$ for all $N\in\N$ and
$\mathbbm{1}^{[N]}\to\mathbbm{1}$, pointwise on $\no$, for $N\to\infty$.
Moreover,
$$
\mathcal{C}_1(\mathbbm{1}-\mathbbm{1}^{[N]})
=\Big(0,0,\dots,0,\frac{1}{N+1},\frac{2}{N+2},\frac{3}{N+3},\dots\Big)
$$
with 0 occurring $N$-times. Accordingly,
$$
\|\mathbbm{1}-\mathbbm{1}^{[N]}\|_{ces_\infty}=\sup_{k\ge1} \frac{k}{k+N}=1,
\quad N\in\N.
$$
So,  $\mathbbm{1}$ does not have a.c.-norm in the B.f.s.\ 
$ces_\infty$, \cite[Proposition I.3.6]{B-S}, and hence, $ces_\infty$ 
does not have a.c.-norm. In particular, $ces_\infty$
is \textit{non-separable}; see Lemma \ref{l-41}.

(ii) Consider the non-negative vector 
$z:=\sum_{k=1}^\infty ke_{k^2}$, which
belongs to $ces_\infty$. Indeed,  for each $n\in\N$ 
let $m^2$ be the largest square integer which is smaller than or equal to $n$.
Then
\begin{align*}
(\mathcal{C}_1z)_n&=\frac{1}{n+1}\sum_{j=0}^n|z_j|
=\frac{1}{n+1}\sum_{j=0}^{m^2}|z_j|
=\frac{1}{n+1}\sum_{k=1}^{m}|z_{k^2}|
\le\frac{1}{m^2}\sum_{k=1}^{m}k
\le1.
\end{align*}
Accordingly, $\mathcal{C}_1|z|\in\ell^\infty$ and so $z\in ces_\infty$. 
Note that $(\mathcal{C}_1z)_0=0$.

Fix any $y=(y_n)_{n=0}^\infty\in\ell^\infty$ and set $M:=\|y\|_\infty$. 
Let $k_0:=2[M]+2$ (where $[x]$ is the integer part of $x$).
We bound 
the norm $\|z-y\|_{ces_\infty}$ from below. Namely,
\begin{align*}
\|z-y\|_{ces_\infty}&=\|\mathcal{C}_1|z-y|\|_{\infty}
=\sup_{n\in\no}\frac{1}{n+1}\sum_{j=0}^n|z_j-y_j|
\\ &  
\ge\frac{1}{4k_0^2+1}\sum_{j=0}^{4k_0^2}|z_j-y_j|
\ge\frac{1}{4k_0^2+1}\sum_{k=1}^{2k_0}|z_{k^2}-y_{k^2}|
\\ & 
=\frac{1}{4k_0^2+1}\sum_{k=1}^{2k_0}|k-y_{k^2}|
\ge\frac{1}{4k_0^2+1}\sum_{k=k_0}^{2k_0}|k-y_{k^2}|
\ge \frac{1}{4k_0^2+1}\sum_{k=k_0}^{2k_0}\big(k-M)
\\ &
=\frac{1}{4k_0^2+1}\left(\frac{(k_0+2k_0)(2k_0-k_0+1)}{2}-(k_0+1)M\right)
\\ &
=\frac{1}{4k_0^2+1}\left(\frac{3k_0(k_0+1)}{2}-(k_0+1)M\right)
\\ &
=\frac{k_0(k_0+1)}{4k_0^2+1}\left(\frac{3}{2}-\frac{M}{k_0}\right)
\\ &
\ge\frac14\left(\frac{3}{2}-\frac12\right)
=\frac14.
\end{align*}
This shows that $z$ does not belong to the closure of
$\ell^\infty$ in $ces_\infty$.
\end{proof}

It was noted prior to  Theorem \ref{t-56}  that all the assumptions
of Theorem \ref{t-35} are satisfied for $X=ces_\infty$ with
the exception that $\ee$ is \textit{not} a basis for $ces_\infty$.
This fact affects the fine spectra of $\ct$ (cf. \eqref{78}) and shows that 
 situation (iv) in the Conjecture can occur.
For the case of $X=\ell^\infty$ we refer to \cite[Theorem 1.4]{S-El}
to observe this feature of the fine spectra of $\ct$.


\begin{theorem}\label{t-77} 
Let  $t\in[0,1)$. The operator $\ct\in\mathcal{L}(ces_\infty)$
satisfies  
$$
\|\ct\|_{ces_\infty\to ces_\infty}= 1,
$$
is compact and has spectrum
\begin{equation}\label{77}
\sigma\big(\ct; ces_\infty\big)=\Lambda\cup\{0\}.
\end{equation}
Moreover, the fine spectra of $\ct$ are given by
\begin{equation}\label{78}
\sigma_{\text{p}}\big(\ct; ces_\infty\big)=\Lambda,
\quad 
\sigma_{\text{c}}\big(\ct; ces_\infty\big)=\emptyset,
\quad
\sigma_{\text{r}}\big(\ct; ces_\infty\big)= \{0\}.
\end{equation}
\end{theorem}


\begin{proof} 
For the continuity and the operator norm of $\ct$ see  Remark \ref{r-57}(i).

By Lemma \ref{l-31} the diagonal operator $D\in \mathcal{L}(ces_\infty)$
is compact with $\|D\|_{ces_\infty\to ces_\infty}=1$ and
by Lemma \ref{l-54} the right-shift operator $S$ satisfies
$\|S\|_{ces_\infty\to ces_\infty}\le1$. So, in view of Lemma  \ref{l-53} if
$x\in ces_\infty\subseteq \cno$, then the series on the right-side of \eqref{25} 
converges absolutely in $ces_\infty$, necessarily to $\ct x\in ces_\infty$;
see Lemma \ref{l-22}. Since $D\in \mathcal{L}(ces_\infty)$ is compact
and the series $\ct=D\sum_\subi t^nS^n$ is absolutely convergent for the
operator norm in $\mathcal{L}(ces_\infty)$, it follows that $\ct$ is compact.

Concerning  the spectrum,
an examination of the proof of Lemma 2.2 in \cite{S-El} shows that
$\sigma_{\text{p}}\big(\ct; ces_\infty\big)\subseteq \Lambda$.
The discussion prior to  Theorem \ref{t-56} shows that $d_1\subseteq ces_\infty$
and hence, via Lemma \ref{l-34}, we can conclude that 
$\Lambda\subseteq \sigma_{\text{p}}\big(\ct; ces_\infty\big)$.
Accordingly, $\sigma_{\text{p}}\big(\ct; ces_\infty\big)= \Lambda$.
This equality, together with the compactness of $\ct\in \mathcal{L}(ces_\infty)$
and the result of F. Riesz concerning the spectrum of compact operators, imply that
\eqref{77} is valid.

Lemma \ref{l-51}(i), for $p=\infty$, shows that $\ct(ces_\infty)\subseteq \ell^\infty$. 
Since the closure of $\ell^\infty$ in $ces_\infty$ is a \textit{proper} subspace
(cf. Proposition \ref{p-76}(ii)) and $\ct$ is injective
(cf.\ Lemma \ref{l-21}(iv)), it follows that 
$0\in \sigma_{\text{r}}\big(\ct; ces_\infty\big)$. That 
$\sigma_{\text{c}}\big(\ct; ces_\infty\big)=\emptyset$ follows immediately.
\end{proof}


\begin{remark}\label{r-78} 
(i) The difference between $\sigma_{\text{c}}\big(\ct; ces_\infty\big)$ and
$\sigma_{\text{r}}\big(\ct; ces_\infty\big)$ to that of
$\sigma_{\text{c}}\big(\ct; ces_p\big)$ and
$\sigma_{\text{r}}\big(\ct; ces_p\big)$, for $p\in\{0\}\cup(1,\infty)$ 
should be noted.

(ii) Since $ces_0\subseteq ces_\infty$ and $ces_0\not\subseteq \ell^\infty$
(cf. Remark  \ref{57}(ii)), it follows that also 
$ces_\infty\not\subseteq \ell^\infty$. According to \eqref{45},
the B.f.s.\  $ces_\infty$ is \textit{not} r.i.

(iii) It follows from \eqref{77} and  Theorem \ref{t-77}(i) that the spectral radius
satisfies
$$
r(\ct;ces_\infty)=\|\ct\|_{ces_\infty\to ces_\infty}=1,\quad t\in[0,1).
$$
\end{remark}


We now turn our attention to $\ct$, for $t\in[0,1)$, acting in another
class of related spaces.
The dual Banach space $(ces_q)'$ of $ces_q$, for  $q\in(1,\infty)$, is rather complicated,
\cite{Ja}. A more tractable \textit{isomorphic} version of $(ces_q)'$, denoted
by $d_p$ with $\frac1p+\frac1q=1$, is presented in Corollary 12.17 of \cite{Be}.
For each  $p\in(1,\infty)$ the space $d_p$ is a reflexive, $p$-convex, 
(complex) Banach lattice
which is continuously contained in $\ell^p$ (for the natural inclusion); see Section 2 of 
\cite{bonet-ricker}.

From the definition of the space $d_1$ (see \eqref{34}) and \eqref{13} it follows 
easily that $e_0\in d_1$ but $\mathcal{C}_1e_0=(\frac{1}{n+1})_\subi\notin d_1$.
Hence, the classical \ces operator $\mathcal{C}_1$ does \textit{not} map $d_1$ into $d_1$.
For $p\in(1,\infty)$ it is known that $\mathcal{C}_1$ does map $d_p$ continuously
into itself and that the spectrum of $\ce_1$ is uncountable; see Proposition
3.2 of \cite{bonet-ricker}. In particular, $\ce_1$ is \textit{not} compact. So, we will
restrict  attention to the operators $\ct$ for $t\in[0,1)$.

For each $p\in[1,\infty)$ define
\begin{equation}\label{61}
d_p:=\Big\{x\in\ell^\infty:\hat x=\Big(\sup_{k\ge n}|x_k|\Big)_\subi\in\ell^p\Big\},
\end{equation}
where $\hat x$  is the least decreasing majorant of $x$; see \eqref{34}.
Then $d_p$ is a Banach lattice for the order from $\cno$ and for the Riesz norm
\begin{equation}\label{62}
\|x\|_{d_p}:=\|\hat x\|_p\, ,\quad x\in d_p .
\end{equation}
Moreover, $d_p$ is \textit{isomorphic} to the dual Banach space $(ces_q)'$, where 
$q\in\{0\}\cup(1,\infty)$ satisfies $\frac1p+\frac1q=1$. For $p\in(1,\infty)$ we refer to
Corollary 12.17 of \cite{Be} and for the fact that $(ces_0)'=d_1$ (with
equality of norms) we refer to Lemma 6.2 of \cite{curbera-ricker4}. Actually,
the isomorphism between $d_p$ and $(ces_q)'$, for $p\in(1,\infty)$, is 
a  Banach lattice
isomorphism; see Proposition 2.4 of \cite{bonet-ricker}. 
It is clear from the definition of $\hat x$ given in \eqref{34} that if 
$x,y\in\ell^\infty$ satisfy $|x|\le|y|$, then $\|x\|_{d_p}\le\|y\|_{d_p}$
(because $\hat x\le \hat y$).
This shows 
that $\|\cdot\|_{d_p}$ is indeed a  Riesz norm in $d_p$, for $p\in[1,\infty)$. Moreover,
if $y\in d_p$ and $x\in\cno$ satisfy $|x|\le|y|$, then also $x\in d_p$.
That is, $d_p$ is a \textit{solid} subspace of $\cno$. 
For each $1\le p<\infty$ the canonical vectors $\ee\subseteq\cno$
belong to $d_p$, form an unconditional basis for the Banach space
$d_p$ and satisfy $\|e_n\|_{d_p}=(n+1)^{1/p}$ for each $n\in\no$;
see Proposition 2.1 of \cite{bonet-ricker} and Section 6 of \cite{curbera-ricker4}.
In particular, $c_{00}\subseteq d_p$.
For further properties of the Banach spaces $d_p$, for $p\in[1,\infty)$
see \cite{Be}, Section 2 of \cite{bonet-ricker} and Section 6 of \cite{curbera-ricker4}.

It follows from \eqref{35} and \eqref{62} that $d_p\subseteq \ell^p$, 
for $p\in[1,\infty)$, with the natural (and proper) inclusion map satisfying 
$\|x\|_p\le \|x\|_{d_p}$, for $x\in d_p$; 
see Proposition 2.7(iii) of \cite{bonet-ricker} for $p\in(1,\infty)$ and note that the
proof also applies for $p=1$. 
Since $\ell^p\subseteq\cno$ continuously, for each $1\le p\le \infty$, it follows that also
$d_p\subseteq \cno$, for  $1\le p< \infty$, with a continuous inclusion.
In particular, $d_p$ is a s.B.l.\ with natural basis for every $p\in[1,\infty)$.
If $x\in\ell^\infty$ is eventually decreasing, then $x\in d_p$ if and only if $x\in \ell^p$ (as
$x\in d_p$ if and only if $\hat x\in \ell^p$).


Given $p\in[1,\infty)$,
define a function norm $\rho_p\colon\mathcal{M}^+\to[0,\infty]$ via
$$
\rho_p(x):=
\begin{cases}
\|\hat x\|_p,\quad x\in \mathcal{M}^+\cap \ell^\infty;
\\
\infty, \quad\text{otherwise}.
\end{cases}
$$
It follows from \eqref{61} and \eqref{62}  that
$$
d_p=X(\rho_p):=\{x\in\cno:\rho_p(|x|)<\infty\}
$$
is a B.f.s.\  on $\no$ 
satisfying $\rho_p(x)=\|x\|_{d_p}$
for $x\in d_p$.
Also the discussion after \eqref{41}  implies that $d_p$ is B.f.s.


Note, for $x\in\cno$ eventually decreasing and  each $p\in[1,\infty)$, that 
\begin{equation*}
x\in d_p \text{ if and only if } x\in\ell^p,
\end{equation*}
and in this case,
\begin{equation*}
\|x\|_p\le \|\hat x\|_{d_p},\quad x\in d_p.
\end{equation*}
%


\begin{lemma}\label{l-62}
For each $p\in[1,\infty)$ the right-shift operator $S\colon d_p\to d_p$  satisfies 
\begin{equation}\label{65}
\|S^m\|_{d_p\to d_p}= (m+1)^{1/p},\quad m\in\no.
\end{equation}
\end{lemma}

\begin{proof}
Fix $p\in[1,\infty)$. Given $m\in\no$, from \eqref{24}  we observe that 
\begin{equation*}
S^mx=(0,0,\dots0,x_0,x_1,x_2,\dots),\quad x\in d_p,
\end{equation*}
with 0 occurring $m$ times before $x_0$. Hence, for $x\in d_p$, we have that
$$
(S^mx)\hat{}=((\hat x)_0,\dots,(\hat x)_0,(\hat x)_1,(\hat x)_2,\dots) 
$$
with $(\hat x)_0$ occurring $(m+1)$ times. It follows that
$$
\big\|(S^mx)\,\hat{}\;\big\|_p^p=(m\cdot ((\hat x)_0)^p)+\|\hat x\|_p^p\le 
(m+1)\|\hat x\|_p^p, 
$$
that is, $\big\|(S^mx)\big\|_{d_p}\le(m+1)^{1/p}\|x\|_{d_p}$. 
Since $x\in d_p$ is arbitrary, we can conclude that 
$\|S^m\|_{d_p\to d_p}\le (m+1)^{1/p}$.

The basis vector  $e_0\in d_p$ satisfies $S^me_0=e_m$ and so 
$(S^me_0)\,\hat {}=\sum_{j=0}^me_j$, for $m\in\no$.
Hence,
$$
\|S^me_0\|_{d_p}^p=\|(S^me_0)\,\hat {}\;\|_{p}^p=
\|\sum_{j=0}^me_j\|_p^p=m+1,\quad m\in\no.
$$
Since $\|e_0\|_{d_p}=1$, we can conclude that
$\|S^m\|_{d_p\to d_p}= (m+1)^{1/p}$, for all $m\in\no$. 
\end{proof}


\begin{corollary}\label{c-63}
For each $p\in[1,\infty)$ and $t\in[0,1)$ the series  
\begin{equation}\label{67}
R_t:=\sum_\subi t^nS^n
\end{equation}
converges absolutely for the operator norm 
in $\mathcal{L}(d_p)$. Moreover, $\|R_0\|_{d_p\to d_p}=1$ and 
\begin{equation*}
(1-t^p)^{-1/p}\le \|R_t\|_{d_p\to d_p}
\le (1-t)^{-1-1/p},\quad t\in(0,1).
\end{equation*}
\end{corollary}

\begin{proof}
Fix $p\in[1,\infty)$ and $t\in[0,1)$. According to \eqref{65} we have that
$$
\sum_\subi \big\|t^nS^n\big\|_{d_p\to d_p}=\sum_\subi t^n\big\|S^n\big\|_{d_p\to d_p}
= \sum_\subi (n+1)^{1/p}t^n<\infty
$$
and so the series \eqref{67} is absolutely convergent in the 
Banach space $(\mathcal{L}(d_p),\|\cdot\|_{d_p\to d_p})$.
An application of H\"older's inequality yields
\begin{align*}
\|R_t\|_{d_p\to d_p}& \le
\sum_\subi (n+1)^{1/p}t^n =\sum_\subi (n+1)^{1/p}t^{n/p}t^{n/q}
\\ & \le \Big(\sum_\subi (n+1)t^{n}\Big)^{1/p}
\Big(\sum_\subi t^{n}\Big)^{1/q}
\\ &=(1-t)^{-2/p}(1-t)^{-1/q}=(1-t)^{-1-(1/p)}.
\end{align*}

Recalling that $S^me_0=e_m$, for $m\in\no$, it follows that
$$
R_te_0=\sum_\subi t^n S^ne_0=\sum_\subi t^n e_n=(1,t,t^2,\dots).
$$
Since $R_te_0$ is a decreasing sequence, we have that 
$(R_te_0)\,\hat{}=R_te_0$ and so
$$
\|R_te_0\|_{d_p}^p=\|(R_te_0)\,\hat{}\,\|_{p}^p=\sum_\subi t^{np}
=(1-t^p)^{-1},
$$
that is, $\|R_te_0\|_{d_p}=(1-t^p)^{-1/p}$. But, $\|e_0\|_{d_p}=1$ and we can conclude 
that
$$
(1-t^p)^{-1/p}=\|R_te_0\|_{d_p}\le \|R_t\|_{d_p\to d_p},\quad t\in(0,1).
$$
For $t=0$ we have that $\|R_0\|_{d_p\to d_p}=1$ as $R_0=I$.
\end{proof}


Each $a\in\cno$ defines a linear \textit{convolution operator} 
$T_a\colon\cno\to\cno$ via
\begin{equation*}
x\mapsto T_ax:=x*a:=\Big(\sum_{j=0}^n x_ja_{n-j}\Big)_\subi,
\quad (x_n)_\subi\in\cno.
\end{equation*}
Note that $x*a=a*x$. Let $p\in[1,\infty)$. If $x*a\in\ell^p$ for each $x\in\ell^p$,
then $a$ is called a \textit{$p$-multiplier} for $\ell^p$. The Closed Graph
Theorem ensures that the restriction of $T_a$ to $\ell^p$ belongs 
to $\mathcal{L}(\ell^p)$. The following result occurs  in \cite{Ni};
see also \cite[p.78]{CMR}.


\begin{proposition}\label{p-64}
Let $a=(a_n)_\subi\in\ell^1$. Then $a$ is a $p$-multiplier for
$\ell^p$ for every $p\in[1,\infty)$, that is, $T_a\in\mathcal{L}(\ell^p)$.
Moreover, $\|a\|_p\le \|T_a\|_{\ell^p\to\ell^p}\le \|a\|_1$.
\end{proposition}


We have an immediate application.

\begin{proposition}\label{p-65}
For  $p\in(1,\infty)$ and $t\in[0,1)$ the  operator
$\ct\in \mathcal{L}(d_p)$. Moreover, $\|\ce_0\|_{d_p\to d_p}=1$
 and, with $\xi:=(1/(n+1))_{n=0}^\infty$, we have that
$$
\|\ct\|_{d_p\to d_p}\le\frac{\|\xi\|_p}{1-t},\quad t\in(0,1).
$$
\end{proposition}

\begin{proof}
Since $\ce_0=D$, the case $t=0$ is covered by Lemma \ref{l-31} with $X=d_p$.

So, let $p\in(1,\infty)$ and $t\in(0,1)$. Fix $x\in d_p$. It follows from \eqref{11} and
\eqref{35} that
\begin{equation}\label{69}
|\ct x|\le \ct|x|=\Big(\frac{1}{n+1}\sum_{k=0}^n t^{n-k}|x_k|\Big)_\subi
\le \Big(\frac{1}{n+1}\sum_{k=0}^n t^{n-k}(\hat x)_k\Big)_\subi.
\end{equation}

For each $n\in\no$ the term $\sum_{k=0}^n t^{n-k}(\hat x)_k$ is the $n$-th coordinate
of $a*\hat x$, where $a:=(t^n)_\subi\in\ell^1$ and $\hat x\in\ell^p$. It follows
from Proposition \ref{p-64} that
$$
\sum_{k=0}^n t^{n-k}(\hat x)_k\le \|a*\hat x\|_p\le \|a\|_1\|\hat x\|_p=(1-t)^{-1}\|x\|_{d_p},
\quad n\in\no.
$$
According to \eqref{69} we can conclude that
$|\ct x|\le (1-t)^{-1}\|x\|_{d_p}\xi$
and hence, that
$\|\ct x\|_{d_p}\le  (1-t)^{-1}\|x\|_{d_p}\|\xi\|_{d_p}$.
Since $\xi$ is a decreasing sequence, it follows that $\hat \xi=\xi$ and so
$\|\xi\|_{d_p}=\|\xi\|_p$. Accordingly,
$$
\|\ct x\|_{d_p}\le  \frac{\|\xi\|_p}{1-t}\|x\|_{d_p},\quad x\in d_p.
$$
The proof is thereby complete.
\end{proof}


\begin{remark}\label{r-66}
Given $t\in[0,1]$ and $p\in(1,\infty)$ an  alternative estimate
for the norm of $\ct\colon d_p\to d_p$ is possible. Indeed, fix $x\in d_p$. Then
\eqref{69}  yields
$$ 
|\ct x|\le \Big(\frac{1}{n+1}\sum_{k=0}^n t^{n-k}|x_k|\Big)_\subi
\le \Big(\frac{1}{n+1}\sum_{k=0}^n |x_k|\Big)_\subi
=\ce_1|x|.
$$
Since $\|\ce_1\|_{d_p\to d_p}=q=p/(p-1)$,
\cite[Proposition 3.2(iv)]{bonet-ricker}, it follows that
$$
\|\ct x\|_{d_p}\le\|\ce_1|x|\|_{d_p}
\le\frac{p}{p-1}\|x\|_{d_p},
$$ 
from which we can conclude that also 
$$
\|\ct\|_{d_p\to d_p}\le\frac{p}{p-1}.
$$
\end{remark}


Since the sequence $\xi\not\in d_1$, the proof of Proposition \ref{p-65}
does  not apply when $p=1$ and $t\in[0,1)$. To include this case we need to
proceed in a different manner.

\begin{proposition}\label{p-67}
Let $p\in[1,\infty)$ and $D\in\mathcal{L}(d_p)$ be the diagonal operator
in \eqref{12}. Then 
\begin{equation}\label{610}
\ct=DR_t=\sum_\subi t^nDS^n,\quad t\in[0,1),
\end{equation}
where $S\in\mathcal{L}(d_p)$ is the right-shift operator  in \eqref{24},
the operator $R_t\in\mathcal{L}(d_p)$ is given by \eqref{67} and the
series in \eqref{610} is absolutely convergent for the operator norm in 
$\mathcal{L}(d_p)$.
\end{proposition}

\begin{proof}
Fix $p\in[1,\infty)$. For $t=0$ the identity \eqref{610} is clear as $\ce_0=D$. So, we can 
assume that $t\in(0,1)$. Since $\|D\|_{\txop}=1$ (cf. Lemma \ref{l-31} with $X=d_p$), it
follows from Corollary \ref{c-63} that the series in \eqref{610} is absolutely
convergent in the Banach space $(\mathcal{L}(d_p),\|\cdot\|_{d_p\to d_p})$
and converges to $DR_t$. 
So, we can apply Proposition \ref{p-32} (see also its proof)
with $X=d_p$ to conclude that  \eqref{610} is valid.
\end{proof}


\begin{remark}\label{r-68}
(i) Note that $p=1$ is included in the statement of Proposition \ref{p-67}
and hence, $\ct\in\mathcal{L}(d_1)$ for all $t\in[0,1)$.

(ii) According to  Corollary \ref{c-63} and 
the fact that $\|D\|_{d_p\to d_p}=1$  it follows from \eqref{610}
that
\begin{equation*}\label{611}
\|\ct\|_{d_p\to d_p} \le (1-t)^{-1-(1/p)},\quad t\in(0,1),
\quad p\in[1,\infty).
\end{equation*}
On the other hand,  Proposition \ref{p-65} yields that also
\begin{equation*}\label{612}
\|\ct\|_{d_p\to d_p} \le \frac{\|\xi\|_p}{(1-t)},\quad t\in(0,1),\quad p\in(1,\infty).
\end{equation*}
The alternative estimate
\begin{equation*}\label{613}
\|\ct\|_{d_p\to d_p} \le \frac{p}{p-1},\quad p\in(1,\infty).
\end{equation*}
is also available (see Remark \ref{r-66}).

(iii) It was shown in the proof of Corollary \ref{c-63}, 
for  $t\in[0,1)$ and $p\in[1,\infty)$, that $R_te_0=\sum_\subi t^ne_n$. 
Via \eqref{25} we can conclude that 
$\ct e_0=DR_te_0=
\big(\frac{t^n}{n+1}\big)_\subi$,
which is a decreasing sequence. Accordingly,
$$
\|\ct e_0\|_{d_p}^p=\|(\ct e_0)\,\hat{}\,\|_p^p=\sum_\subi\frac{t^{np}}{(n+1)^p}.
$$
Since $\|e_0\|_{d_p}=1$, we can conclude for each $1<p<\infty$ that
\begin{equation}\label{614}
1<\Big(\sum_\subi\frac{t^{np}}{(n+1)^p}\Big)^{1/p}\le \|\ct\|_{d_p\to d_p},
\quad t\in(0,1).
\end{equation}
\end{remark}


Since all the assumptions of Theorem \ref{t-35}, for $X=d_p$\,, have 
been verified, we can deduce the following result.

\begin{theorem}\label{t-69} 
Let $t\in[0,1)$ and  $p\in[1,\infty)$. 
The operator $\ct\in\mathcal{L}(d_p)$ satisfies 
$$
\|\ct\|_{d_p\to d_p}\le (1-t)^{-1-(1/p)},
$$ 
is compact and has spectrum
$$
\sigma(\ct;d_p)=\Lambda\cup\{0\}.
$$
Moreover, the fine spectra of $\ct$ are given by
$$
\sigma_{\text{p}}(\ct;d_p)=\Lambda,\quad
\sigma_{\text{c}}(\ct;d_p)=\{0\},\quad 
\sigma_{\text{r}}(\ct;d_p)=\emptyset.
$$
\end{theorem}


\begin{remark}\label{r-610} 
(i) According to \eqref{45}, the B.f.s.\ $X=d_1$ is \textit{not} r.i.\ because it
does not contain $\ell^1$. Indeed, the element
$$
x:=\Big(1,0,\frac12,0,0,0,\frac{1}{2^2},\dots,\frac{1}{2^{n-1}},0,\dots,0,\frac{1}{2^{n}},0\Big),
$$
where $0$ occurs $2^n-1$ times between $\frac{1}{2^{n-1}}$ and $\frac{1}{2^{n}}$ 
surely belongs to $\ell^1$. On the other hand, its least decreasing majorant
$$
\hat x:=\Big(1,\frac12,\frac12,\frac{1}{2^2},\frac{1}{2^2},\frac{1}{2^2},\frac{1}{2^2},
\dots,\frac{1}{2^{n-1}},\frac{1}{2^{n}},\dots,\frac{1}{2^{n}},\frac{1}{2^{n+1}},\dots\Big),
$$
where $\frac{1}{2^{n}}$ occurs $2^n$ times, does not belong to $\ell^1$, 
that is, $x\notin d_1$. For $X=d_p$ with $p\in(1,\infty)$ it follows from 
Remark 2.8(i) in \cite{bonet-ricker} that also $\ell^1\not\subseteq d_p$ and so
$d_p$ also  \textit{fails} to be a r.i.\ space.

(ii) According to Theorem \ref{t-69} the spectral radius $r(\ct;d_p)=1$ for
$1\le p<\infty$ and $t\in[0,1)$. Combined with \eqref{614} we can conclude
that
$$
r(\ct;d_p)<\|\ct\|_{d_p\to d_p},\quad 1<p<\infty,\; t\in(0,1).
$$
\end{remark}


\section{The operators $\ct$ acting in weighted $\ell^p$  and $c_0$ spaces}
\label{S8}


Let $w=(w(n))_\subi$ be a strictly positive, decreasing sequence satisfying 
$w(n)\downarrow0$. Define the \textit{weighted space}
$$
\ell^p(w):=\Big\{x\in\cno: \|x\|_{p,w}:=\Big(\sum_\subi |x_n|^pw(n)\Big)^{1/p}<\infty\Big\},
\quad p\in[1,\infty),
$$
equipped with the norm $\|\cdot\|_{p,w}$. Being linearly
and isometrically isomorphic to $\ell^p$, via the map
$x\mapsto (w(n)^{1/p}x_n)_\subi$, for $x\in \ell^p(w)$, 
it follows that $\ell^p(w)$ is a Banach space and  $\ee\subseteq \ell^p(w)$
is a basis. It is routine to verify, for the  order
induced via $\cno$, that $\ell^p(w)$ is a 
Banach lattice for the Riesz norm $\|\cdot\|_{p,w}$ and that $\ell^p(w)$ is solid in $\cno$.
Observe, for $x\in\ell^p(w)$, that
$|x_n|\le \|x\|_{p,w}/w(n)^{1/p}$, for $n\in\no$.
So, for a given $n\in\no$, the Riesz seminorm in \eqref{23} satisfies
$r_n(x)=\max_{0\le k\le n}|x_k|\le M_n \|x\|_{p,w}$, for $x\in\ell^p(w)$,
where $M_n:=\max_{\,0\le k\le n}w(k)^{-1/p}=w(n)^{-1/p}$ (as $w$ is decreasing).
It follows that the natural inclusion map $\ell^p(w)\subseteq\cno$ is continuous for
each $p\in[1,\infty)$. 
So, for $p\in[1,\infty)$ and for every strictly positive weight $w$ with $w(n)\downarrow0$,
the weighted space $\ell^p(w)$
is a s.B.l.\ with natural basis.
Actually, for each $p\in[1,\infty)$ and each weight $w$ as above, 
$\ell^p(w)$ is also a B.f.s.\  on $\no$
(as it satisfies the Fatou property).

Since $d_1\subseteq \ell^1$, to verify that $d_1\subseteq\ell^p(w)$ it suffices
to show that $\ell^1\subseteq \ell^p(w)$. So, fix $x\in\ell^1$. Since 
$w$ is decreasing, it follows that
$$
\|x\|_{p,w}^p=\sum_\subi|x_n|^pw(n)\le w(0) \sum_\subi|x_n|^p
=w(0)\|x\|_p^p\le w(0)\|x\|_1^p.
$$
Hence, $x\in\ell^p(w)$ and $\|x\|_{p,w}\le w(0)^{1/p}\|x\|_1$.
Accordingly, $\ell^1\subseteq \ell^p(w)$.

To establish that $\ell^p(w)$ is translation invariant fix $x\in\ell^p(w)$.
For the right-shift operator $S$ (cf.  \eqref{24}) it follows that
$\|S\|_{\ell^p(w)\to\ell^p(w)}\le1$ because
\begin{equation}\label{71}
\|Sx\|_{p,w}=\big(\sum_\subi|x_n|^pw(n+1)\Big)^{1/p}\le
\big(\sum_\subi|x_n|^pw(n)\Big)^{1/p}=\|x\|_{p,w},
\end{equation}
for all $x\in \ell^p(w)$. This inequality yields 
\begin{equation}\label{72}
\sum_\subi t^n\|S^n\|_{\ell^p(w)\to\ell^p(w)}\le 
\sum_\subi t^n=(1-t)^{-1}<\infty,\quad t\in[0,1).
\end{equation}
Since
$$
\frac{\|Se_n\|_{p,w}}{\|e_n\|_{p,w}}=\frac{w(n+1)^{1/p}}{w(n)^{1/p}},
\quad n\in\no,
$$
it follows that the weight $w$ necessarily satisfies
$$
\sup_{n\in\no}\frac{w(n+1)^{1/p}}{w(n)^{1/p}}\le \|S\|_{\ell^p(w)\to\ell^p(w)}\le1.
$$

Lemma \ref{l-31} implies that the diagonal operator $D$ in \eqref{12}
belongs to $\mathcal{L}(\ell^p(w))$, satisfies 
$\|D\|_{\ell^p(w)\to\ell^p(w)}=1$ and is compact.
Proposition \ref{p-32} and its proof 
(together with \eqref{32} and \eqref{72})
show that  $\|\ct\|_{\ell^p(w)\to\ell^p(w)}\le(1-t)^{-1}$ and 
that $\ct$ is a compact operator. Then,  
for $X=\ell^p(w)$, Theorem \ref{t-35} implies  the following result.


\begin{theorem}\label{t-71} 
Let $p\in[1,\infty)$ and $t\in[0,1)$. Let $w=(w(n))_\subi$ be any strictly positive
weight on $\no$ satisfying $w(n)\downarrow0$. 
The   operator $\ct\in\mathcal{L}(\ell^p(w))$ satisfies
\begin{equation}\label{73}
\|\ct\|_{\ell^p(w)\to\ell^p(w)}\le (1-t)^{-1},
\end{equation}
is compact and has spectrum
$$
\sigma(\ct;\ell^p(w))=\Lambda\cup\{0\}.
$$
Moreover, the fine spectra of $\ct$ are given by
$$
\sigma_{\text{p}}(\ct;\ell^p(w))=\Lambda,\quad
\sigma_{\text{c}}(\ct;\ell^p(w))=\{0\},\quad 
\sigma_{\text{r}}(\ct;\ell^p(w))=\emptyset.
$$
\end{theorem}


\begin{remark}\label{r-72}
(i) The right-side of  \eqref{73} does not involve $p\in[1,\infty)$
or  $w$.

(ii) The B.f.s.\  $\ell^p(w)$  \textit{fails} to be r.i.\ relative to $\mu$
for each $p\in[1,\infty)$ and each
strictly positive weight $w=(w(n))_\subi$ satisfying $w(n)\downarrow0$.
Indeed, since $w(n)\downarrow0$ with $w(n)>0$ for all $n\in\no$, there exists
$m\in\no$ satisfying $w(m+1)<w(m)$. Consider the canonical vector
$e_m\in\ell^p(w)$. We have noted earlier that $e_m$ and $Se_m$
are equimeasurable. Hence, if $\ell^p(w)$ is r.i., then 
$\|Se_m\|_{p,w}=\|e_m\|_{p,w}$. But, it follows from \eqref{71} that
$\|Se_m\|_{p,w}=w(m+1)^{1/p}<w(m)^{1/p}=\|e_m\|_{p,w}$.

(iii) For  $\ce_1\in\mathcal{L}(\ell^p(w))$, we refer to
\cite{ABR1} when $p\in(1,\infty)$ and to \cite{ABR2} when $p=1$.
\end{remark}


Still under the assumption that   $w=(w(n))_\subi$ is a strictly positive weight 
with $w(n)\downarrow0$, define the weighted space 
$$
c_0(w):=\Big\{x\in\cno: \lim_{n\to\infty} w(n)|x_n|=0\Big\},
$$
equipped with the norm $\|x\|_{0,w}:=\sup_nw(n)|x_n|$ for $x\in c_0(w)$. Being linearly
and isometrically isomorphic to $c_0$, via the map
$x\mapsto (w(n)x_n)_\subi$, for $x\in c_0(w)$,
it follows that $c_0(w)$ is a Banach space and that $\ee\subseteq c_0(w)$
is a basis. Moreover, for the coordinate-wise order, $c_0(w)$ is a 
Banach lattice for the Riesz norm $\|\cdot\|_{0,w}$ 
and $c_0(w)$ is solid in $\cno$.
The inequalities
$|x_n|\le \|x\|_{0,w}/w(n)$, for $n\in\no$,
valid for each $x\in c_0(w)$, imply that the natural inclusion  map
$c_0(w)\subseteq \cno$ is continuous. So,  $c_0(w)$ is a
s.B.l.\ with natural basis.

Since $w\in\ell^\infty$, it is clear that $\ell^1\subseteq c_0(w)$
and hence,   also $d_1\subseteq c_0(w)$. To see that $X:=c_0(w)$ is 
\textit{translation invariant} observe  that
$\|S\|_{c_0(w)\to c_0(w)}\le1$ because
\begin{equation}\label{74}
\|Sx\|_{0,w}=\sup_{n\in\no}w(n+1)|x_n|\le \sup_{n\in\no}w(n)|x_n|=\|x\|_{0,w},
\quad x\in c_0(w).
\end{equation}
It follows from \eqref{74}, together with 
$$
\frac{\|Se_n\|_{0,w}}{\|e_n\|_{0,w}}=\frac{w(n+1)}{w(n)},
\quad n\in\no,
$$
that the weight $w$ necessarily satisfies
$$
\sup_{n\in\no}\frac{w(n+1)}{w(n)}\le \|S\|_{c_0(w)\to\ c_0(w)}\le1.
$$
Note that 
\begin{equation}\label{75}
\sum_\subi t^n\|S^n\|_{c_0(w)\to\ c_0(w)}\le \sum_\subi t^n=(1-t)^{-1},\quad t\in[0,1).
\end{equation}
Lemma \ref{l-31} implies  the diagonal operator $D$ in \eqref{12} belongs 
to $\mathcal{L}(c_0(w))$, satisfies $\|D\|_{c_0(w)\to c_0(w)}=1$ and is compact.
Proposition \ref{p-32} and its proof 
(together with \eqref{33} and\eqref{75})
show that $\ct\in\mathcal{L}(c_0(w))$
with $\|\ct\|_{c_0(w)\to c_0(w)}\le(1-t)^{-1}$ and that $\ct$ is compact. 
Applying Theorem \ref{t-35}, for $X=c_0(w)$,  yields the following result.


\begin{theorem}\label{t-73} 
Let $ t\in[0,1)$ and $w=(w_n)_\subi$ be any strictly positive
weight on $\no$ satisfying $w(n)\downarrow0$.
The   operator $\ct\in\mathcal{L}(c_0(w))$ satisfies
\begin{equation}\label{76}
\|\ct\|_{c_0(w)\to\ c_0(w)}\le (1-t)^{-1},
\end{equation}
is compact and has spectrum
$$
\sigma(\ct;c_0(w))=\Lambda\cup\{0\}.
$$
Moreover, the fine spectra of $\ct$ are given by
$$
\sigma_{\text{p}}(\ct;c_0(w))=\Lambda,\quad
\sigma_{\text{c}}(\ct;c_0(w))=\{0\},\quad 
\sigma_{\text{r}}(\ct;c_0(w))=\emptyset.
$$
\end{theorem}


\begin{remark}\label{r-74}
(i) The right-side of  \eqref{76} does not involve the weight $w$.

(ii) The space $c_0(w)$ \textit{fails} to be a B.f.s.\ 
on $\no$ as it does not have the Fatou property.

(iii) For  the relevant properties of $\ce_1\in\mathcal{L}(c_0(w))$, we refer to
\cite{ABR3}.
\end{remark}


\section{The operators $\ct$ acting in Bachelis spaces}
\label{S9}


Let $\T=\{z\in\C:|z|=1\}$ denote the boundary (in $\C$) of $\D=\{z\in\C:|z|<1\}$.
For each $1<p<\infty$ define
$$
N^p(\Z):=\Big\{x\in\C^\Z: |x|\le \hat f \text{ for some } f\in L^p(\T)\Big\},
$$
where $\hat f$ denotes the Fourier transform of $f$. Then 
$\np$ is a reflexive Banach space relative to the norm
\begin{equation}\label{671}
\|x\|_\np:=\inf\Big\{\|f\|_p: |x|\le \hat f,\; f\in L^p(\T)\Big\},
\end{equation}
where $\|\cdot\|_p$ is the norm in $L^p(\T)$, that is,
$\|f\|_p=\big(\frac1{2\pi}\int_{-\pi}^\pi|f(e^{it})|^p\,dt\big)^{1/p}$.
The  Banach spaces $\np$, $1<p<\infty$, were introduced
by Bachelis, \cite{Ba}. Since $L^p(\T)$ is reflexive, the infimum
in \eqref{671} is achieved. Let $P_p\in\mathcal{L}(\np)$ denote the projection
$$
P_p x:=(x_n)_\subi,\quad x=(x_n)_{n\in\Z}\in\np.
$$
Its range $N^p:=P_p(\np)$ is a complemented subspace of $\np$ identifiable with
$$
N^p\simeq \Big\{x\in\cno: |x|\le \hat f,\; f\in H^p\Big\},
$$
after recalling that every function $f$ in the Hardy space
$H^p=H^p(\D)$ is holomorphic in $\D$, its boundary function (again
denoted by $f$) belongs to $L^p(\T)$ and satisfies
$\hat f(n)=0$ for every $n\in\Z$ with $n<0$.
Moreover, if $f\in L^p(\T)$ and we define 
$F(z):=\sum_{n=0}^\infty \hat f(n)z^n$, then $F\in H^p$ and 
every $F\in H^p$ is so obtained. For these facts see
\cite[Chapter 17]{Ru}, for example.

The canonical vectors $\mathcal{E}$ form a (unconditional) basis
in $N^p$, for $1<p<\infty$. Moreover, the dual Banach space of $N^p$ can be
identified isometrically with $N^q$, where $\frac1p+\frac1q=1$, via the pairing
$\langle x,u\rangle=\sum_\subi x_nu_n$ for $x\in N^p$ and $u\in N^q$,
which arises from the corresponding dual pairing between
$N^p(\Z)$ and $N^q(\Z)$, \cite[Theorem 2]{Ba}. Each space $N^p$ is a reflexive 
Banach lattice for the coordinate-wise order inherited from $\cno$ (with
respect to which it  is solid) and the Riesz norm 
\begin{equation}\label{672}
\|x\|_{N^p}:=\inf\Big\{\|f\|_p: |x|\le \hat f,\; f\in H^p\Big\},
\quad x\in N^p.
\end{equation}
Moreover, $N^p\subseteq c_0$ continuously and hence, also the natural
inclusion $N^p\subseteq\cno$ is continuous. Accordingly, 
$N^p$ is a s.B.l.\  with natural basis.

Given $x\in\ell^1$ the function $f(z)=\sum_\subi|x_n|z^n$ is holomorphic
in $\D$ and continuous in $\overline{\D}$. Hence, $f\in H^p$ and satisfies
$|x|\le \hat f$. Since $|f|\le \|x\|_1\chi_\T$ and $\|\chi_\T\|_p=1$, it follows
(cf. \eqref{672}) that
$$
\|x\|_{N^p}\le \|f\|_p\le \|x\|_1.
$$
So, $\ell^1\subseteq N^p$ with a continuous inclusion and hence, also
$d_1\subseteq N^p$ continuously.

\begin{lemma}\label{l-71}
Let $p\in(1,\infty)$.
\begin{itemize}
\item[(i)] The basis $\mathcal{E}\subseteq N^p$ is normalized, that is, $\|e_n\|_{N^p}=1$ for all $n\in\no$.
\item[(ii)] The shift operator $S\colon N^p\to N^p$ satisfies $\|S^m\|_{N^p\to N^p}= 1$
for all $m\in\no$.
\end{itemize}
\end{lemma}

\begin{proof}
(i) Fix $p\in[2,\infty)$. Given $n\in\no$ the function $f_n(z)=z^n$, for $z\in\D$, belongs
to $H^p$ and satisfies $\|f_n\|_p=1$ with $|e_n|\le \hat {f_n}$ (as $\hat{f_n}=\chi_{\{n\}}$).
It follows from \eqref{672} that $\|e_n\|_{N^p}\le1$. On the other hand,
consider any function $f\in H^p$ satisfying $|e_n|\le \hat {f}$, in which case
$1\le \hat f(n)$. Since the inclusion $L^p(\T)\subseteq L^2(\T)$ has operator norm 1,
it follows from Plancherel's Theorem that
$$
\|f\|_p\ge \|f\|_2=\Big(\sum_{k=0}^\infty |\hat f(k)|^2\Big)^{1/2}\ge f(n)\ge1.
$$
Then \eqref{672} yields that $\|e_n\|_{N^p}\ge1$. This proves that $\|e_n\|_{N^p}=1$.

Now let $p\in(1,2)$. Fix $n\in\no$. The same argument as for $p\in[2,\infty)$
yields that $\|e_n\|_{N^p}\le1$. Since $N^p$ is isometrically isomorphic to
$N^q$ (with $\frac1p+\frac1q=1$), it follows that
$$
\|e_n\|_{N^p}=\sup\big\{|\langle e_n,u\rangle|: u\in N^q,\; \|u\|_{N^q}=1\big\}.
$$
Choose $u:=e_n\in N^q$. Since $q\in[2,\infty)$ we know that $\|u\|_{N^q}=1$
and hence, $\|e_n\|_{N^p}\ge |\langle e_n,u\rangle|=1$. So, again $\|e_n\|_{N^p}=1$.

(ii) Fix $p\in (1,\infty)$ and let $x\in N^p$. Choose a function $f\in H^p$ satisfying 
$\|x\|_{N^p}=\|f\|_p$ and $|x|\le \hat f$, that is, $|x_n|\le \hat f(n)$ for $n\in\no$.
The function $g(z)=zf(z)$, for $z\in\D$, belongs to $H^p$ and satisfies 
$\hat g(n)=\hat f(n-1)$ for $n\in\Z$. In particular, $\hat g(n)=0$ for all $n\in\Z$
with $n\le0$. Moreover, $\|g\|_p=\|f\|_p$.
Observe  that
$$
|Sx|=S|x|=(0,|x_0|,|x_1|,\dots)\le(0,\hat f(0),\hat f(1),\dots)
=(\hat g(0),\hat g(1),\dots),
$$
that is, $|Sx|\le \hat g$, Via \eqref{672} we can conclude that
$\|Sx\|_{N^p}\le \|g\|_p =\|f\|_p=\|x\|_{N^p}$.
Since $x\in N^p$ is arbitrary, we can deduce that $\|S\|_{N^p\to N^p}\le 1$.
Given $m\in\N$, it follows that $\|S^m\|_{N^p\to N^p}\le 1$.
On the other hand $S^me_0=e_m$ and so, by part (i),
we see that $\|S^m\|_{N^p\to N^p}\ge \|S^me_0\|_{N^p}/\|e_0\|_{N^p}=1$.
Hence, $\|S^m\|_{N^p\to N^p}= 1$.
\end{proof}

According to Lemma \ref{l-71} and the discussion prior to that result, 
we can conclude that each $N^p$ ($1<p<\infty$) is a translation invariant s.B.l.\  
with natural basis.

For $X=N^p$ ($1<p<\infty$), Lemma \ref{l-71}(ii) implies that the series
in \eqref{32} satisfies
$$
\sum_\subi\|t^nS^n\|_{N^p\to N^p}\le \sum_\subi t^n=\frac{1}{1-t}<\infty,\quad t\in[0,1),
$$
and so this series is absolutely convergent for the operator norm in 
$\mathcal{L}(N^p)$, with limit $R_t$. Since $\|D\|_{N^p\to N^p}=1$ (cf. Lemma \ref{l-31})
and $\ct=D R_t$, we can conclude that $\|\ct\|_{N^p\to N^p}\le(1-t)^{-1}$
for $t\in[0,1)$. Hence, for $X=N^p$, all the assumptions of Theorem \ref{t-35}
are satisfied and so we have established the following result.


\begin{theorem}\label{t-671} 
Let $t\in[0,1)$ and $p\in(1,\infty)$.
The   operator $\ct\in\mathcal{L}(N^p)$ satisfies
$$
\|\ct\|_{N^p\to N^p}\le(1-t)^{-1},
$$
is compact and has spectrum
$$
\sigma(\ct;N^p)=\Lambda\cup\{0\}.
$$
Moreover, the fine spectra of $\ct$ are given by
$$
\sigma_{\text{p}}(\ct;N^p)=\Lambda,\quad
\sigma_{\text{c}}(\ct;N^p)=\{0\},\quad 
\sigma_{\text{r}}(\ct;N^p)=\emptyset.
$$
\end{theorem}

\begin{remark}\label{r-672}
Concerning $t=1$, it follows from \cite[Theorem 2.5]{curbera-ricker5}
that
$$
\|\ce_1\|_{N^p\to N^p}=p,\quad p\in[2,\infty) \text{ and } 
p\le \|\ce_1\|_{N^p\to N^p}\le 2,\quad p\in(1,2).
$$
Furthermore, $\sigma_{\text{p}}(\ce_1;N^p)=\emptyset$ for all $1<p<\infty$. Moreover,
$$
\sigma(\ce_1;N^p)=\Big\{z\in\C:\big|z-\frac{p}{2}\big|\le\frac{p}{2}\Big\},
\quad p\in[2,\infty),
$$
and
$$
\Big\{z\in\C:\big|z-\frac{p}{2}\big|\le\frac{p}{2}\Big\}
\subseteq \sigma(\ce_1;N^p),\quad p\in(1,2).
$$
In particular, $\ce_1 \in\mathcal{L}(N^p)$ is not a compact
operator for $1<p<\infty$.
\end{remark}


\section{Concluding remarks}
\label{S77}


We have seen that Theorem \ref{t-35} is rather versatile and applies to a large
and varied class of s.B.l.'s $X$,
thereby providing strong evidence for the validity of the Conjecture in Section 1.
The central idea behind Theorem \ref{t-35} is the factorization of the operator
\begin{equation}\label{n71}
C_t=DR_t,\quad t\in[0,1),
\end{equation}
with $D\in\mathcal{L}(X)$ being compact and $R_t\in \mathcal{L}(X)$
having a representation
\begin{equation}\label{n72}
R_t=\sum_\subi t^nS^n,
\end{equation}
where the series is convergent for the strong operator topology in 
$\mathcal{L}(X)$ and $S\in \mathcal{L}(X)$ is the right-shift operator;
see Proposition \ref{p-32}. This approach requires $X$ to be translation invariant;
see \eqref{n72}. We now present an example of a s.B.l.\ $X$ with
natural basis which contains $d_1$ but \textit{fails} to be translation
invariant. Accordingly, Theorem \ref{t-35} is not applicable.
Nevertheless, each operator $\ct\in \mathcal{L}(X)$, for $t\in[0,1)$,
turns out to be compact and has spectrum and fine spectra as stated 
in the Conjecture.

For $x\in\cno$ let $\xe:=(x_{2n})_\subi\in\cno$ consist
of those coordinates of $x$ having an even subscript and
$\xo:=(x_{2n+1})_\subi\in\cno$ consist
of those coordinates of $x$ having an odd subscript.
Fix any pair $1<p<q<\infty$ and equip the linear space
\begin{equation*}
X_{p,q}:=\Big\{x\in\cno: \xe\in\ell^p,\; \xo\in\ell^q\Big\}
\end{equation*}
with the Riesz norm
\begin{equation*}
\|x\|_{p,q}:=\|\xe\|_p+ \|\xo\|_q,\quad x\in X_{p,q},
\end{equation*}
in which case $X_{p,q}$ is a s.B.l.\ for the coordinate-wise
order from $\cno$. Observe that $\mathcal{E}\subseteq X_{p,q}$
satisfies $\|e_n\|_{p,q}=1$ for all $n\in\no$.
It is routine to verify that $\mathcal{E}$ is a basis for $X_{p,q}$.
Since $\ell^q\subseteq\cno$ continuously, Lemma \ref{l-n71} below implies that 
 $X_{p,q}$ is a s.B.l.\ with natural basis.
Moreover, $d_1\subseteq \ell^1$ and $\ell^1$ is contained 
in both $\ell^p$ and $\ell^q$, from which it follows that $d_1\subseteq X_{p,q}$.

Let $z\in\ell^q\setminus\ell^p$ and define $x=(0,z_0,0,z_1,0,z_2,0\dots)$,
in which case $\xe=0\in\ell^p$ and $\xo=z\in\ell^q$, that is, $x\in X_{p,q}$.
Note that
$$
Sx=(0,0, z_0,0,z_1,0,z_2,0,\dots).
$$
Since $(Sx)^{(\text{e})}=(0,z_0,z_1,z_2\dots)\notin\ell^p$, it is clear 
that $Sx\notin X_{p,q}$. Accordingly, $X_{p,q}$ is \textit{not} 
translation invariant. In particular, the B.f.s.\ $X_{p,q}$ cannot be r.i.

\begin{lemma}\label{l-n71}
For each pair $1<p<q<\infty$ it is the case that 
\begin{equation}\label{n73}
\ell^p\subsetneq X_{p,q}\subsetneq \ell^q
\end{equation}
with both inclusions continuous and having operator norm
at most 2.
\end{lemma}

The proof follows routinely from the definition of $\|\cdot\|_{p,q}$ and is omitted.

\begin{lemma}\label{l-n72}
Let $t\in[0,1)$ and $1<p<q<\infty$. The operator
$\ct\in\mathcal{L}(X_{p,q})$ satisfies 
\begin{equation}\label{n74}
\|\ct\|_{X_{p,q}\to X_{p,q}}\le 2\|\xi\|_p (1-t)^{-1},
\end{equation}
where $\xi:=(1/(n+1))_\subi$.
\end{lemma}

\begin{proof}
Let $x\in X_{p,q}$. Since both
$\xe,\xo\in\ell^q$, it follows that
\begin{equation}\label{n74-75}
\|x\|_q=(\|\xe\|_q^q+\|\xo\|_q^q)^{1/q}
\le (\|\xe\|_p^q+\|\xo\|_q^q)^{1/q}
\le \|\xe\|_p+\|\xo\|_q=\|x\|_{p,q}.
\end{equation}
From \eqref{11} it is clear that
\begin{equation}\label{n75}
|\ct x|\le \ct|x|=\Big(\frac{1}{n+1}\sum_{k=0}^n t^{n-k}|x_k|\Big)_\subi.
\end{equation}
Fix $n\in\no$. The term $\sum_{k=0}^n t^{n-k}|x_k|=(a*|x|)_n$
is the $n$-th coordinate of $a*|x|$, where
$a=(t^n)_\subi\in\ell^1$ and $|x|\in\ell^q$ (cf.\ \eqref{n73}).
It follows from \eqref{n74-75} and Proposition \ref{p-65}  that
$$
\sum_{k=0}^n t^{n-k}|x_k|\le \|a*|x|\|_q\le \|a\|_1\,\|x\|_q
\le \|a\|_1\,\|x\|_{p,q}=(1-t)^{-1}\|x\|_{p,q}.
$$
According to \eqref{n75} we can conclude that
$$
|\ct x|\le (1-t)^{-1}\|x\|_{p,q}\; \xi,
$$ 
and hence, that
$$
\|\ct x\|_{p,q}\le (1-t)^{-1}\|x\|_{p,q} \,\|\xi\|_{p,q}.
$$ 
Since $\xi\in\ell^p\subseteq X_{p,q}$ and
$\|\xi\|_{p,q}=\|\xi^{\text{(e)}}\|_p+\|\xi^{\text{(o)}}\|_q
\le \|\xi^{\text{(e)}}\|_p+\|\xi^{\text{(o)}}\|_p
\le 2\|\xi\|_p$,
it follows  that
$$
\|\ct x\|_{p,q}\le 2(1-t)^{-1} \|\xi\|_{p}\,\|x\|_{p,q}.
$$ 
This  inequality implies \eqref{n74}.
\end{proof}

Since $X_{p,q}$ fails to be translation invariant, it is clear that
$\ct\in\mathcal{L}(X_{p,q})$ \textit{cannot} have a factorization
of the form \eqref{n71}. Nevertheless, the following result shows that
$\ct$ is still compact.

\begin{lemma}\label{l-n73}
For each $t\in[0,1)$ and $1<p<q<\infty$ the operator
$\ct\in\mathcal{L}(X_{p,q})$ is compact.
\end{lemma}

\begin{proof}
Fix $t\in[0,1)$ and $1<p<q<\infty$. Let $N\in\N$. Define the finite-rank operator
$\ct^{[N]}\in\mathcal{L}(X_{p,q})$ by
$$
\ct^{[N]} x:= \sum_{n=0}^N\frac{1}{n+1} (a*x)_n e_n=
\sum_{n=0}^N\frac{1}{n+1} \Big(\sum_{k=0}^n t^{n-k} x_k\Big)e_n,
\quad x\in X_{p,q},
$$
with $a:=(t^n)_\subi\in\ell^1$. Setting $\xi^{[N]}:=\sum_{m=N+1}^\infty \xi_me_m\in
\ell^p\subseteq X_{p,q}$ (with $\xi$ as in Lemma \ref{l-n72}) we
have, for each $x\in X_{p,q}$, that
\begin{align*}
\Big|\big(\ct-\ct^{[N]}\big)x\Big|=
\Big|\big(\xi_n^{[N]}(a*x)_n\big)_\subi\Big|
\le \Big(\xi_n^{[N]}(a*|x|)_n\Big)_\subi
\le (1-t)^{-1} \|x\|_{p,q} \;\xi^{[N]};
\end{align*}
see the proof of Lemma \ref{l-n72}. It follows that
$$
\Big\|\big(\ct-\ct^{[N]}\big)x\Big\|_{p,q}\le  2(1-t)^{-1}\|\xi^{[N]}\|_p\;\|x\|_{p,q}.
$$
Since $x\in X_{p,q}$ is arbitrary, the previous inequality yields
$$
\Big\|\big(\ct-\ct^{[N]}\big)\Big\|_{X_{p,q}\to X_{p,q}}\le  2(1-t)^{-1}\|\xi^{[N]}\|_p,
\quad N\in\N.
$$
But, $\|\xi^{[N]}\|_p\to0$ for $N\to\infty$ (as $\xi\in\ell^p$) and so we can conclude
that $\ct$ is compact.
\end{proof}

An examination of the proof of Theorem \ref{t-35} reveals that it can be
adapted to establish the following result.

\begin{theorem}\label{t-n74}
Let $X$ be a  s.B.l.\
with  natural basis which contains $d_1$.
Suppose, for some $t\in[0,1)$, that  $\ct\colon X\to X$  exists and  is a compact operator. 
Then its spectrum
$$
\sigma(\ct;X)=\Lambda\cup\{0\}. 
$$
Moreover, the fine spectra of $\ct$ are given by
$$
\sigma_{\text{p}}(\ct;X)=\Lambda,\quad 
\sigma_{\text{c}}(\ct;X)=\{0\},\quad
\sigma_{\text{r}}(\ct;X)=\emptyset.
$$
\end{theorem}

We end with two applications of Theorem \ref{t-n74}.

\begin{corollary}\label{c-nn75}
Let $X$ be a s.B.l.\ with natural basis which contains $d_1$ and
such that both $X$ and $X'$ have a.c.-norm. 
Suppose that $\ce_{t_0}\colon X\to X$ exists and is a 
compact operator for some $t_0\in(0,1)$. Then
$\ct\in\mathcal{L}(X)$ is compact for every $t\in[0,t_0]$ and its
spectrum 
$$
\sigma(\ct;X)=\Lambda\cup\{0\}. 
$$
Moreover, the fine spectra of $\ct$ are given by
$$
\sigma_{\text{p}}(\ct;X)=\Lambda,\quad 
\sigma_{\text{c}}(\ct;X)=\{0\},\quad
\sigma_{\text{r}}(\ct;X)=\emptyset.
$$
\end{corollary}

\begin{proof}
Since $C_0=D$, the case $t=0$ is clear.

In view of Theorem \ref{t-n74} it suffices to show, for any fixed $t\in(0,t_0)$, that
$\ct\colon X\to X$ exists and is compact. Concerning the existence,
let $x\in X$. By Lemma \ref{l-21}(ii) we have $0\le|\ct x|\le \ce_{t_0}|x|$ with 
$\ce_{t_0}|x|\in X$. Since $X$ is solid, also $\ct x\in X$. Accordingly, 
$\ct\colon X\to X$ exists. Since $X\subseteq\cno$ continuously and
$\ct\in\mathcal{L}(\cno)$, a closed graph argument shows that 
$\ct\in\mathcal{L}(X)$.

Let $X_{\R}:=X\cap \R^{\no}$. Then $X_\R$ is a real Banach lattice
and $X=X_\R\oplus i X_\R$ is the complexification
of $X_\R$ (cf.~\cite[II Section 11]{Sch}, for example). Let
$\widetilde \ct\colon X_\R\to X_\R$ be the restriction of $\ct$ to $X_\R$
which, according to \eqref{11}, takes its values in $X_\R$.
Again by \eqref{11} it follows that $\widetilde \ct\ge0$ is a positive operator in $X_\R$ and hence, $\ct\ge0$ is a positive operator in the complex Banach
lattice $X$, \cite[p.135]{Sch}. The norms in $X$ and $X'$, both being
a.c., imply that the same is true for the restrictions of these norms to
$X_\R$ and $(X')_\R=(X_\R)'$ respectively, after recalling that 
$X'=(X_\R)'\oplus i (X_\R)'$. Moreover, it is clear from \eqref{11} and
Lemma \ref{l-21}(ii) that $0\le \widetilde \ct\le \widetilde{\ce_{t_0}}$ as positive
operators in $X_\R$. Since a positive operator in $X$ is compact if and only of its restriction to $X_\R$ is compact (as can be seen by splitting sequences in $X$ into
real and imaginary parts), it follows that 
$\widetilde\ce_{t_0}\in\mathcal{L}(X_\R)$ is compact. Hence, by
the Dodds-Fremlin theorem, \cite[Theorem 3.7.13]{MN}, we can conclude that
$\widetilde\ct\in\mathcal{L}(X_\R)$ is compact from which the
compactness of $\ct\in\mathcal{L}(X)$ follows.
\end{proof}

In view of Lemma \ref{l-n72}, Lemma \ref{l-n73} and the discussion prior
to Lemma \ref{l-n71} we have the following consequence of Theorem
\ref{t-n74}. Recall that $X_{p,q}$ is \textit{not} translation invariant.

\begin{corollary}\label{c-n75}
Let $t\in[0,1)$ and $1<p<q<\infty$. The  operator
$\ct\in\mathcal{L}(X_{p,q})$ satisfies 
$$
\|\ct\|_{X_{p,q}\to X_{p,q}} \le \frac{2\|\xi\|_p}{1-t},
$$
is  compact and has  spectrum
$$
\sigma(\ct;X_{p,q})=\Lambda\cup\{0\}. 
$$
Moreover, the fine spectra of $\ct$ are given by
$$
\sigma_{\text{p}}(\ct;X_{p,q})=\Lambda,\quad 
\sigma_{\text{c}}(\ct;X_{p,q})=\{0\},\quad
\sigma_{\text{r}}(\ct;X_{p,q})=\emptyset.
$$
\end{corollary}



\end{document}